\documentclass{article}

\usepackage{amsmath,amssymb, amsthm}
\usepackage{physics}
\usepackage[]{comment}
\usepackage{adjustbox}
\usepackage{subcaption, siunitx, booktabs}
\usepackage{tikz}
\usepackage{mathtools}
\usepackage{parskip}
\usepackage{pifont}
\usepackage{graphicx}
\usepackage{hyperref}
\usepackage{authblk}

\DeclareMathOperator*{\argmin}{arg\,min}

\title{Higher-order meshless schemes for hyperbolic equations}
\author[1, 2]{Klaas Willems\thanks{klaas.willems@kuleuven.be, \url{https://orcid.org/0009-0004-8527-4621}}}
\author[2]{Axel Klar\thanks{axel.klar@rptu.de, \url{https://orcid.org/0000-0003-1679-9438}}}
\author[1]{Giovanni Samaey\thanks{giovanni.samaey@kuleuven.be, \url{https://orcid.org/0000-0003-1679-9438}}}
\affil[1]{Department of Computer Science, KU Leuven, Leuven, Belgium}
\affil[2]{Faculty of Mathematics, RPTU Kaiserslautern-Landau, Kaiserslautern, Germany}

\providecommand{\keywords}[1]
{
	\small	
	\textbf{Keywords. } #1
}

\providecommand{\classification}[1]
{
	\small	
	\textbf{MSC subject classifications. } #1
}

\begin{document}

\newtheorem{remark}{Remark}
\newcommand{\cmark}{\ding{51}}%
\newcommand{\xmark}{\ding{55}}%

\maketitle

\begin{abstract}
	We discuss the order, efficiency, stability and positivity of several meshless schemes for linear scalar hyperbolic equations. Meshless schemes are Generalised Finite Difference Methods (GFDMs) for arbitrary irregular grids in which there is no connectivity between the grid points. We propose a new MUSCL-like meshless scheme that uses a central stencil, with which we can achieve arbitrarily high orders, and compare it to existing meshless upwind schemes and meshless WENO schemes. The stability of the newly proposed scheme is guaranteed by an upwind reconstruction to the midpoints of the stencil. The new meshless MUSCL scheme is also efficient due to the reuse of the GFDM solution in the reconstruction. We combine the new MUSCL scheme with a Multi-dimensional Optimal Order Detection (MOOD) procedure to avoid spurious oscillations at discontinuities. In one spatial dimension, our fourth order MUSCL scheme outperforms existing WENO and upwind schemes in terms of stability and accuracy. In two spatial dimensions, our MUSCL scheme achieves similar accuracy to an existing WENO scheme but is significantly more stable.
\end{abstract}
 \\
\classification{65M06, 35L04}

\section{Introduction}
Finite element and finite volume methods have been used with great success in many fields of study in both academic and industrial circles. These methods are however not without limitations. Generating three-dimensional meshes to be used in these methods can be a time consuming task. In addition, due to the reliance on a mesh, finite volume methods and finite element methods are not well suited for problems with moving boundaries. One strategy for dealing with these issues is to remesh \cite{belytschko_meshless_1996}. This is however costly (again), and requires the projection of quantities between consecutive meshes \cite{nguyen_meshless_2008}. 

Several approaches have been developed over the years to overcome the difficulties associated with mesh dependence. One such approach is based on Radial Basis Functions (RBFs). RBF methods approximate numerical data through local or global interpolations using radially symmetric basis functions. They were first introduced for solving partial differential equations in \cite{kansa_multiquadricsscattered_1990, kansa_multiquadricsscattered_1990-1}, and have since been applied to a wide range of problems; see \cite{fornberg_solving_2015} and references therein. An alternative approach are so-called Generalised Finite Difference Methods (GFDMs) \cite{jensen_finite_1972, perrone_general_1975, liszka_finite_1980}. The origins of GFDMs can be traced back to smoothed particle hydrodynamics \cite{monaghan_introduction_1988, lucy_numerical_1977}, which was initially developed for modelling astrophysical phenomena, and was later used in fluid dynamics \cite{gomez-gesteira_state---art_2010}. GFDMs only require local connectivity information: for each grid point a set of nearby points is required, and not complete volumes or elements. As a result, these methods work well in combination with Lagrangian or arbitrary Lagrangian-Eulerian (ALE) schemes in which the grid movement is determined by the equations \cite{tiwari_meshfree_2022}. Moreover, due to the absence on any constraints on the grid, generation of the grids is cheap. \\

The fundamental idea of GFDMs is to approximate spatial derivatives in a partial differential equation using a polynomial reconstruction through the function values at the surrounding points. This approximation is based on a weighted least squares method \cite{shepard_two-dimensional_1968, lancaster_surfaces_1981, levin_approximation_1998} known as Moving Least Squares (MLS). The main differences between the plethora of GFDMs proposed in literature lies in the selection of the stencil, the polynomial reconstruction, the weight function used in the least squares problem, and any additional constraints that are added to the least squares problem. Based on these choices, we can distinguish several types of GFDMs. In `classical' GFDMs \cite{benito_influence_2001, benito_solving_2007, gavete_solving_2017, prieto_application_2011, urena_solving_2019}, the polynomial reconstruction is applied at each grid point, and is based on a Taylor expansion. The least squares problem provides a relation between the spatial derivatives at a point, and the neighbouring function values. This expression can then be used during time integration. In the Finite Pointset Method (FPM) \cite{kuhnert_grid_2001, tiwari_finite_2003, suchde_meshfree_2018, resendiz-flores_two-dimensional_2015}, also Taylor expansions are used for the polynomial reconstruction. The main difference with respect to classical GFDMs, is the inclusion of the physical equations in the weighted least squares problem using a Lagrange multiplier. Finally, in the Finite Point Method \cite{boroomand_simple_2005, boroomand_generalized_2009}, the polynomial reconstruction is computed on subcells of the domain, not for each grid point separately. See \cite{clain_stencil_2024}, for an excellent overview of GFDMs. \\

The majority of the research on GFDMs has been for parabolic and elliptic partial differential equations \cite{kuhnert_grid_2001, benito_influence_2001}, for problems in elasticity \cite{liszka_finite_1980}, heat transfer \cite{resendiz-flores_two-dimensional_2015, gavete_solving_2017}, and fluid dynamics \cite{suchde_meshfree_2018}. There exist significantly less results in which GFDMs are used for hyperbolic conservation laws. Ramesh and Deshpande combined an upwinding procedure with a classical GFDM to obtain a first order stable scheme \cite{ramesh_least_2001}. This scheme was later extended to second order \cite{ghosh_least_1995}. One major issue with these schemes is that in two dimensions, they are not positive, leading to spurious oscillations at discontinuities. Chandrashekar obtained a first order positive scheme by adding additional diffusion to the upwind flux \cite{chan_positivity_2021}. The more recent on publications GFDMs for hyperbolic conservation laws focus on the issue of conservation: GFDMs lack formal conservation. Several strategies have been proposed to solve this issue. Chiu et al. \cite{chiu_conservative_2011, kwan-yu_chiu_conservative_2012} add additional constraints to the least squares problem to enforce conservation at the discrete level. In doing so, they lose the local nature of the differential operators, and significantly increase the computational cost of the method. Pratik et al. \cite{suchde_flux_2017} also enforce conservation, but only a local grid. Although their method thus lacks formal conservation, they report significantly reduced errors. Similarly, Huh et al. \cite{huh_new_2018} enforce a property referred to as `geometric conservation', which can prevent the creation or reduction of mass.

Besides the issue of conservation, existing methods are at most of order two, are not positivity-preserving, and their stability properties are not well understood. In the this work, we develop a stable higher-order positive scheme for linear hyperbolic equations on irregular grids. We achieve stability at higher-orders by using a higher-order upwind reconstruction at the midpoints of the stencil, similarly to classical MUSCL schemes using in Finite Volume Methods (FVMs). This upwind reconstruction is then used in a classical GFDM. The reconstruction uses the same central stencil as the GFDM, as a result, the method is very efficient. To obtain a positivity-preserving method, we use the Multi-dimensional Optimal Order Detection (MOOD) procedure \cite{clain_high-order_2011, diot_improved_2012, diot_multidimensional_2013}. In MOOD, the proposed solution at the next time step is checked against a local discrete maximum property (DMP). In case the DMP is not satisfied, the order of the spatial discretisation is reduced to order one, ensuring positivity. We also note that the MUSCL-like scheme can be combined with modern approaches to make GFDMs (approximately) conservative. We compare our approach to a second-order Weighted Essentially Non-Oscillating (WENO) MLS method developed in \cite{tiwari_meshfree_2022} based on \cite{avesani_new_2014}. This WENO MLS method computes an MLS approximation of the gradient using an upwind and central stencil. These results are combined using non-linear weights based on so-called oscillation indicators, which assign a larger weight to smoother stencils. This method has been shown to yield good results in both smooth regions and regions that contain shocks.

This paper is organised as follows. In section \ref{section:GFDMIntro}, the basic least squares procedure of a GFDM is given. In section \ref{section:FirstOrder}, we discuss positivity-preserving first order meshless schemes in one and two spatial dimensions. These are used as a stable `fallback' method for the meshless MUSCL scheme in the MOOD procedure. In section \ref{section:HigherOrder}, we discuss the newly developed meshless MUSCL scheme, and give a meshless WENO scheme. In section \ref{section:Numerics}, we illustrate the properties of the meshless schemes in several numerical experiments. 

\section{Model problem and least squares approximation}
\label{section:GFDMIntro}
As a model problem, we consider the linear advection equation \begin{align}
    \pdv{u}{t} + a \cdot \nabla u = 0 \label{eq:linearAdvection}, 
\end{align} with $a \in \mathbb{R}^d$ a constant vector (\(d = 1 \text{ or } 2\)). An analysis of the schemes for the linear advection equation is sufficient for use with kinetic equations, and is necessary before applying the schemes to non-linear hyperbolic equations. We solve equation \eqref{eq:linearAdvection} on an irregular grid with grid points numbered from $i = 1\dots N$. We denote the position of a grid point as $\vec{x_i}$ and the numerical solution of \eqref{eq:linearAdvection} at that grid point as $u_i$. In addition, for each grid point $i$, we define the set of neighbours $\mathcal{C}_i$. Two points are neighbours if the distance between them is smaller than a threshold value $h_{max}$. Thus, unlike in a uniform grid, stencils do not have a fixed topology, and some grid points may have more neighbours than others. For the following, we do not use any connectivity of the grid, except the neighbourhood relation mentioned above.

\subsection{Moving least squares method}
\label{section:leastSquares}

The method of choice for evaluating gradients on irregular grids is the so-called Moving Least Squares (MLS) method \cite{levin_approximation_1998, lancaster_surfaces_1981}. We summarize the method in the 2D case below. For a more detailed discussion on MLS, see \cite{seibold_m-matrices_2006, suchde_conservation_2018, tey_moving_2021}. MLS relies on a Taylor expansion around the central point $u(x_i, y_i) = u_i$ \begin{align}
    u_j = u_i + \Delta x_{ij} \pdv{u_i}{x} + \Delta y_{ij} \pdv{u_i}{y} + \mathcal{O}(\Delta x_{ij}^2) + \mathcal{O}(\Delta y_{ij}^2), \; \; j \in \mathcal{C}_i, \label{eq:taylorExpansion}
\end{align} where $\Delta x_{ij} = x_j - x_i$ and $\Delta y_{ij} = y_j - y_i$. An approximation of the spatial derivatives is then obtained by minimizing the $L^2$ norm of the interpolation error with respect to a weight function \begin{align}
	\left( \pdv{\Tilde{u}_i}{x}, \pdv{\Tilde{u}_i}{y} \right) = \argmin_{\pdv{\Tilde{u}_i}{x}, \pdv{\Tilde{u}_i}{y}} \quad \sum_{j \in \mathcal{C}_i} w_{ij} \left( u_j -  u_i - \Delta x_{ij} \pdv{\Tilde{u}_i}{x} - \Delta y_{ij} \pdv{\Tilde{u}_i}{y} \right)^2. \label{eq:minLeastSquares}
\end{align} In this text, we use the weight function \begin{align}
	w_{ij} = w(\vec{x}_i, \vec{x}_j) = \exp( -\alpha || \vec{x}_i - \vec{x}_j ||^2 ) \label{eq:WeightFunction}
\end{align} with \(\alpha\) a dimension-dependent parameter, although other weight functions are also possible (see \cite{seibold_m-matrices_2006}). The solution to the minimization problem leads to the following matrix equation \begin{align}
    \begin{bmatrix}
        \sum_{j \in \mathcal{C}_i} w_{ij} \Delta x_{ij}^2 & \sum_{j \in \mathcal{C}_i} w_{ij} \Delta x_{ij} \Delta y_{ij} \\
        \sum_{j \in \mathcal{C}_i} w_{ij} \Delta x_{ij} \Delta y_{ij} & \sum_{j \in \mathcal{C}_i} w_{ij} \Delta y_{ij}^2 \\
    \end{bmatrix} \begin{bmatrix}
        \pdv{\Tilde{u}_i}{x} \\ \pdv{\Tilde{u}_i}{y}
    \end{bmatrix} = \begin{bmatrix}
        \sum_{j \in \mathcal{C}_i} w_{ij} \Delta x_{ij} (u_j - u_i) \\ \sum_{j \in \mathcal{C}_i} w_{ij} \Delta y_{ij} (u_j - u_i)
    \end{bmatrix}.
\end{align} The solution to these equations can be written explicitly as \begin{align}
    \pdv{u_i}{x} \approx \pdv{\Tilde{u}_i}{x} = \sum_{j \in \mathcal{C}_i} \alpha_{ij} (u_j - u_i), \; \pdv{u_i}{y} \approx \pdv{\Tilde{u}_i}{y} = \sum_{j \in \mathcal{C}_i} \beta_{ij} (u_j - u_i), \label{eq:2DDerivatives}
\end{align} where the coefficients $\alpha_{ij}$ and $\beta_{ij}$ are given by \begin{align}
    \alpha_{ij} &= \frac{ \left( \sum_{k \in \mathcal{C}_i} w_{ik} \Delta y_{ik}^2 \right) w_{ij} \Delta x_{ij} -  \left( \sum_{k \in \mathcal{C}_i} w_{ij} \Delta y_{ik} \Delta x_{ik} \right) w_{ij} \Delta y_{ij}}{ \left( \sum_{j \in \mathcal{C}_i} w_{ij} \Delta x_{ij}^2 \right) \left( \sum_{j \in \mathcal{C}_i} w_{ij} \Delta y_{ij}^2 \right) - \left( \sum_{j \in \mathcal{C}_i} w_{ij} \Delta x_{ij} \Delta y_{ij} \right)^2 } \label{eq:2Dalfaij}\\
    \beta_{ij} &= \frac{ \left( \sum_{k \in \mathcal{C}_i} w_{ik} \Delta x_{ik}^2 \right) w_{ij} \Delta y_{ij} -  \left( \sum_{k \in \mathcal{C}_i} w_{ij} \Delta y_{ik} \Delta x_{ik} \right) w_{ij} \Delta x_{ij}}{ \left( \sum_{j \in \mathcal{C}_i} w_{ij} \Delta x_{ij}^2 \right) \left( \sum_{j \in \mathcal{C}_i} w_{ij} \Delta y_{ij}^2 \right) - \left( \sum_{j \in \mathcal{C}_i} w_{ij} \Delta x_{ij} \Delta y_{ij} \right)^2 }. \label{eq:2Dbetaij}
\end{align} The above procedure can be extended to arbitrary order by including more terms in the Taylor expansion \eqref{eq:taylorExpansion}. The extension to 1D and 3D is also straightforward. The only requirement for the above procedure to be used is that the least square problem is solvable. This is the case if there are sufficient points in the stencil (depending on the dimension and order), and the points do not lie on a straight line. 

\subsection{Note on conservation}
\label{section:conservationRemark}
It is well known that the use of the procedure outlined in section \ref{section:leastSquares} for hyperbolic conservation laws does not yield conservative schemes. This is due to two reasons. Firstly, classical meshfree methods lack flux conservation, i.e., the mass moving from one grid point $i$ to grid point $j$ is not equal to the negative of the mass moving from grid point $j$ to grid point $i$. Secondly, meshless schemes don't use control volumes, which allows the method to create or remove mass. The first property can be dealt with by solving a least squares problem analogous to \eqref{eq:minLeastSquares} but for all points at once, with additional equality constraints that enforce flux conservation. The second property can be dealt with by adding so-called `geometric-conservation' constraints to the least squares problem \eqref{eq:minLeastSquares}. In \cite{kwan-yu_chiu_conservative_2012}, the authors developed a scheme that is both flux and geometric conserving, and therefore a conservative scheme in the classical sense. However, due to the increased computational cost, this scheme is practically not usable any more in three-dimensional domains and in the case of moving boundaries, which undermines the merit of meshless schemes. Instead, many authors have chosen to only enforce geometric conservation \cite{katz_meshless_2009, huh_new_2018}. In \cite{suchde_flux_2017}, a local flux conserving method was introduced using a background mesh. 

This work intends to develop arbitrary high-order stable meshless schemes. For this reason, we limit ourselves to the classical generalised finite difference methods as outlined in section \ref{section:leastSquares}. As a result, our method is not conservative. However, the method can be made conservative by using any of the techniques developed in \cite{kwan-yu_chiu_conservative_2012, suchde_flux_2017, katz_meshless_2009, huh_new_2018}.   

\section{First order schemes}
\label{section:FirstOrder}
The meshless MUSCL method developed in this paper uses a Multi-dimensional Optimal Order Detection (MOOD) procedure to enforce positivity of the solution. Such a procedure relies on a positivity-preserving first-order scheme to fall back to in regions where a higher-order method yields spurious oscillations. In this section we discuss a positivity-preserving first-order method in one and two-spatial dimensions.

\subsection{Formulation in one space dimension}
\label{section:FirstOrder1D}
We generalise the classical upwind scheme to general unstructured grids. To this end, we define the upwind stencil \( U_i \) as the set of points that lie `upwind' of \(x_i\), and within a distance \(h_{max}\) of \(x_i\). Using the least squares procedure from section \ref{section:leastSquares} with the stencil \(U_i\) we obtain the following first order semi-discretized meshless scheme \begin{align}
	\pdv{u_i}{t} &= - a \frac{\sum_{k \in U_i} w_{ik} \Delta x_{ik} \left( u_k - u_i \right)}{\sum_{k \in U_i} w_{ik} \Delta x_{ik}^2}, 
\end{align} By simple manipulation of the expression, it can be seen that under forward Euler time integration, this scheme is positive (and thus adheres to the maximum principle) under the following CFL condition \begin{align}
	1 + a \Delta t \frac{\sum_{k \in U_i}^m w_{ik} \Delta x_{ik}}{\sum_{k \in U_i}^m w_{ik} \Delta x_{ik}^2} \geq 0. \label{eq:1DCFLCondition}
\end{align} 

\subsection{Formulation in two space dimensions}
\label{section:FirstOrder2D}

\begin{figure}
	\begin{subfigure}{.32\linewidth}
		\centering
		\resizebox{\linewidth}{!}{
			\begin{tikzpicture}
				\draw[->] (-5, 0) -- (5, 0) node[right] {$x$};
				\draw[->] (0, -5) -- (0, 5) node[above] {$y$};
				
				\fill[blue] (1, 2) circle (2pt);     
				\fill[blue] (-2, 3) circle (2pt);    
				\fill[blue] (-3, -2) circle (2pt);   
				\fill[blue] (3, -1.5) circle (2pt);  
				\fill[blue] (2, 1.5) circle (2pt);     
				\fill[blue] (-2.5, 0.5) circle (2pt);    
				\fill[blue] (-0.5, -1) circle (2pt);   
				\fill[blue] (0.5, -1) circle (2pt);  
				
				\draw[->, thick, red] (0,0) -- (3,4) node[anchor=south west] {$a$};
				\draw[thick] (-5, 3.75) -- (5, -3.75) node[anchor=north west] {};
				
				\draw[dashed] (0, 0) -- (2, 1.5);
				
				\draw[thick] (1,0) arc[start angle=0, end angle=36.87, radius=1cm];
				\node at (1.3,0.2) {$\theta_{ij}$};
				
				\fill[green!20, opacity=0.5] (-5, 3.75) -- (-5,-5) -- (5,-5) -- (5,-3.75) -- (0, 0) -- cycle;
			\end{tikzpicture}
		}
		
		\caption{}
		\label{fig:2DUpwindTrue}
	\end{subfigure}%
	\begin{subfigure}{.32\linewidth}
		\centering
		\resizebox{\linewidth}{!}{
			\begin{tikzpicture}
				\draw[->] (-5, 0) -- (5, 0) node[right] {$x$};
				\draw[->] (0, -5) -- (0, 5) node[above] {$y$};
				
				\fill[blue] (1, 2) circle (2pt);     
				\fill[blue] (-2, 3) circle (2pt);    
				\fill[blue] (-3, -2) circle (2pt);   
				\fill[blue] (3, -1.5) circle (2pt);  
				\fill[blue] (2, 1.5) circle (2pt);     
				\fill[blue] (-2.5, 0.5) circle (2pt);    
				\fill[blue] (-0.5, -1) circle (2pt);   
				\fill[blue] (0.5, -1) circle (2pt);  
				
				\draw[->, thick, red] (0,0) -- (3,4) node[anchor=south west] {$a$};
				\draw[thick] (-5, 3.75) -- (5, -3.75) node[anchor=north west] {};
				
				\draw[dashed] (0, 0) -- (2, 1.5);
				
				\draw[thick] (1,0) arc[start angle=0, end angle=36.87, radius=1cm];
				\node at (1.3,0.2) {$\theta_{ij}$};
				
				\fill[green!20, opacity=0.5] (-5, 0) -- (-5,-5) -- (0,-5) -- (0, 0) -- cycle;
			\end{tikzpicture}
		}
		\caption{}
		\label{fig:2DUpwindQuadrant}
	\end{subfigure}
	\begin{subfigure}{0.32\linewidth}
		\centering
		\resizebox{\linewidth}{!}{
			\begin{tikzpicture}
				\draw[->] (-5, 0) -- (5, 0) node[right] {$x$};
				\draw[->] (0, -5) -- (0, 5) node[above] {$y$};
				
				\fill[blue] (1, 2) circle (2pt);     
				\fill[blue] (-2, 3) circle (2pt);    
				\fill[blue] (-3, -2) circle (2pt);   
				\fill[blue] (3, -1.5) circle (2pt);  
				\fill[blue] (2, 1.5) circle (2pt);     
				\fill[blue] (-2.5, 0.5) circle (2pt);    
				\fill[blue] (-0.5, -1) circle (2pt);   
				\fill[blue] (0.5, -1) circle (2pt);  
				
				\draw[->, thick, red] (0,0) -- (3,4) node[anchor=south west] {$a$};
				\draw[thick] (-5, 3.75) -- (5, -3.75) node[anchor=north west] {};
				
				\draw[dashed] (0, 0) -- (2, 1.5);
				
				\draw[thick] (1,0) arc[start angle=0, end angle=36.87, radius=1cm];
				\node at (1.3,0.2) {$\theta_{ij}$};
				
				\fill[green!20, opacity=0.5] (-5, 5) -- (-5,-5) -- (0,-5) -- (0, 5) -- cycle;
				\fill[orange!20, opacity=0.5] (-5, 0) -- (-5,-5) -- (5,-5) -- (5, 0) -- cycle;
			\end{tikzpicture}
		}
		\caption{}
		\label{fig:2DUpwindSplit}
	\end{subfigure}
	\caption{Possible upwind stencils in two spatial dimensions. The central point of the stencil is the point at the origin. The points in blue are the neighbours \(C_i\). The vector \(a\) is the velocity in the linear advection equation. The shaded regions indicate possible upwind stencils. The angle \(\theta_{ij}\) is the angle between the positive x-axis and the line connecting the centerpoint and the neighbour \(j\).}
	\label{fig:2DUpwind}
\end{figure}
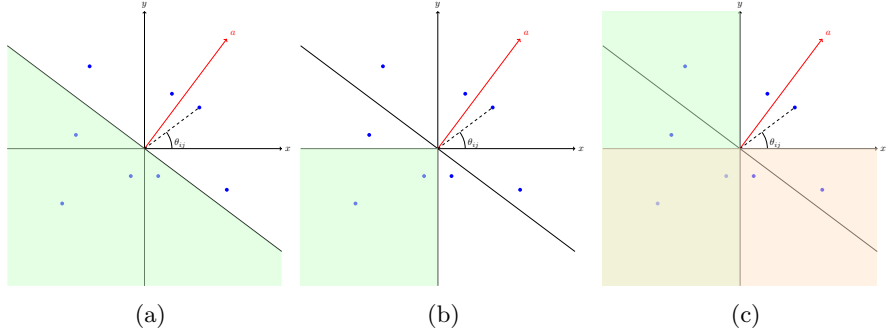

In one dimension, if the velocity is positive, it is obvious the upwind direction is to the left. This idea can be extended to 2D in several ways. One could choose the upwind stencil as all the points `behind' the central point (see \ref{fig:2DUpwindTrue}), or one could restrict the stencil to the upwind quadrant (see \ref{fig:2DUpwindQuadrant}). Another option would be to compute the derivatives in the \(x\) and \(y\)-direction using separate stencils, for example using a left and right half for the x-direction and a top and bottom half for the y-direction (see \ref{fig:2DUpwindSplit}). Using the expressions \eqref{eq:2Dalfaij} and \eqref{eq:2Dbetaij} it is possible to show that none of these upwind methods can yield a positive time integration scheme. Instead, we use a scheme based on a central stencil, to which artificial diffusion is added, that was developed in \cite{chandrashekar_positive_2004}. For this scheme, it is possible to prove positivity under a CFL condition.

We summarise the scheme from \cite{chandrashekar_positive_2004} here. We write the divergence \eqref{eq:linearAdvection} \begin{align}
	\div{a u_i} = 2 a_x \sum_{j \in C_i} \alpha_{ij} (u_{ij} - u_i) + 2 a_y \sum_{j \in C_i} \beta_{ij} (u_{ij} - u_i), \label{eq:PraveenDivergence1}
\end{align} where $\alpha_{ij}$ and $\beta_{ij}$ are given by \eqref{eq:2DDerivatives}. The value $u_{ij}$ is the numerical solution at the midpoint $\frac{x_i + x_j}{2}$. Note that the value \(u_{ij}\) will in fact never be explicitly computed, but only the flux \( a u_{ij} \) at the point. Below, we will give a suitable definition for this flux. The factor two in \eqref{eq:PraveenDivergence1} naturally arises due to the definition of the midpoint. 

Let $\theta_{ij}$ be the angle between the line connecting the grid point $i$ and grid point $j$ and the positive x-axis. The orthogonal vectors then $\hat{n}_{ij} = \left( \cos \theta_{ij}, \sin \theta_{ij} \right)^T$ and $\hat{s}_{ij} = \left( - \sin \theta_{ij}, \cos \theta_{ij} \right)^T$ define a rotated right-handed coordinate system. We rewrite the velocity vector in this coordinate system \begin{align}
	\begin{bmatrix}
		a_x \\ a_y
	\end{bmatrix} = \begin{bmatrix}
		\hat{n}_{ij} & \hat{s}_{ij}
	\end{bmatrix} \begin{bmatrix}
		\left( \vec{a} \cdot \hat{n}_{ij} \right) \\
		\left( \vec{a} \cdot \hat{s}_{ij} \right)
	\end{bmatrix}, \label{eq:PraveenVelRotation}
\end{align} and define \begin{align}
	\begin{bmatrix}
		\Bar{\alpha}_{ij} \\ \Bar{\beta}_{ij}
	\end{bmatrix} = \begin{bmatrix}
		\hat{n}_{ij} & \hat{s}_{ij}
	\end{bmatrix}^T \begin{bmatrix}
		\alpha_{ij} \\ \beta_{ij}
	\end{bmatrix} \label{eq:PraveenCoeffRotation}.
\end{align} Using equations \eqref{eq:PraveenVelRotation} and \eqref{eq:PraveenCoeffRotation}, the divergence in \eqref{eq:PraveenDivergence1} can be written as \begin{align}
	\div{a u_i} = 2 \sum_{j \in C_i} \Bar{\alpha}_{ij} \left( \vec{a} \cdot \hat{n}_{ij} \right) (u_{ij} - u_i) + 2 \sum_{j \in C_i} \Bar{\beta}_{ij} \left( \vec{a} \cdot \hat{s}_{ij} \right) (u_{ij} - u_i). \label{eq:PraveenDivergence2}
\end{align} The first and second term in \eqref{eq:PraveenDivergence2} represent the flux along $\hat{n}_{ij}$ and $\hat{s}_{ij}$. Now, the fluxes at the midpoint \(u_{ij}\) along $\hat{n}_{ij}$ and $\hat{s}_{ij}$ are defined as \begin{align}
	\left( \vec{a} \cdot \hat{n}_{ij} \right) u_{ij} &= \frac{ \vec{a} \cdot \hat{n}_{ij} + |\vec{a} \cdot \hat{n}_{ij}| }{2} u_i + \frac{ \vec{a} \cdot \hat{n}_{ij} - |\vec{a} \cdot \hat{n}_{ij}| }{2} u_j \label{eq:nFlux}\\
	\left( \vec{a} \cdot \hat{s}_{ij} \right) u_{ij} &= \left( \vec{a} \cdot \hat{s}_{ij} \right) \frac{u_i + u_j}{2} - \text{sign}(\Bar{\beta}_{ij})\frac{|\vec{a} \cdot \hat{s}_{ij}|}{2}. \label{eq:sFlux}
\end{align} Equation \eqref{eq:nFlux} is the flux at the midpoint flowing orthogonal to an imaginary cell boundary separating nodes $i$ and $j$. A simple upwind flux is chosen. Equation \eqref{eq:sFlux} is the flux at the midpoint flowing tangential to the boundary. In a finite volume scheme, due to Gauss' theorem, this flux is zero. Here, due to the lack of cells, or any connectivity in general, this flux is not zero. The first term in \eqref{eq:sFlux} is a simple average of the fluxes from the nodes. The second term in \eqref{eq:sFlux} is a diffusive term that forces the scheme to be positive. These choices of flux yield the following semi-discretised scheme
\begin{align}
	\dv{u_i}{t} &= - \sum_{j \in C_i} \left[ \Bar{\alpha}_{ij} \left( \vec{a} \cdot \hat{n}_{ij} - |\vec{a} \cdot \hat{n}_{ij}| \right) + \left( \Bar{\beta}_{ij} \left( \vec{a} \cdot \hat{s}_{ij} \right) - |\Bar{\beta}_{ij} \left( \vec{a} \cdot \hat{s}_{ij} \right)| \right) \right] (u_j - u_i),
\end{align} which in combination with a forward Euler time integration routine, can be shown to be positive under a CFL condition \cite{chandrashekar_positive_2004}. This algorithm was extended to three-dimensional space in \cite{mariappan_high-performance_2024}.

\section{Higher order methods}
\label{section:HigherOrder}
In section \ref{section:MUSCLMOOD}, we introduce a strategy with which arbitrary higher-order stable meshless methods can be constructed. The scheme combines a MUSCL-like reconstruction with MOOD \cite{clain_high-order_2011}. MOOD is an explicit \textit{a posteriori} check of physical constraints such as positivity of density and numerical constraints such as the discrete maximum property. MOOD breaks away from classical \textit{a priori} limiting schemes and WENO schemes. In section \ref{section:WENO}, we summarise the WENO-like method introduced in \cite{avesani_new_2014} for comparison with the new MUSCL-like method. Both methods are given in one spatial dimension, and then extended to 2-D.

\subsection{Meshless MUSCL method with MOOD}
\label{section:MUSCLMOOD}

\subsubsection{Formulation in one space dimension}
\label{section:MUSCLMOOD1D}
In principle, one could obtain a higher-order meshless scheme starting from the scheme in section \ref{section:FirstOrder}, by using upwind stencils and higher-order Taylor expansions \eqref{eq:taylorExpansion}. However, this has several disadvantages. Firstly, since only points from one direction are considered, larger neighbourhoods are required such that the least squares problem remains solvable. Secondly, in higher dimensions, selecting the upwind stencil can be cumbersome, and it is not obvious how narrow to take the upwind direction. Intuitively, narrow upwind directions are more stable, yet require even larger neighbourhoods. Thirdly, in kinetic problems, when one has to solve the PDE for multiple velocities, one would have to solve a least squares problem for each velocity. For these reasons, we have opted to use central stencils instead. However, as we show in section \ref{section:Stability}, a naïve central GFDM yields an unconditionally unstable scheme. We therefore borrow an idea from classical MUSCL schemes, namely to reconstruct the solution at a cell boundary. We give the second order scheme here, and after, discuss the extension to higher order. Consider a third order MLS approximation in one spatial dimension analogous to section \ref{section:leastSquares}. This procedure yields a set of coefficients \(\kappa_{ij}\) and \(\lambda_{ij}\) with which we can compute approximations of the first order and second order spatial derivatives. \begin{align}
	\pdv{u_i}{x} \approx \pdv{\Tilde{u}_i}{x} = \sum_{j \in \mathcal{C}_i} \kappa_{ij} (u_j - u_i), \; \; \pdv[2]{u_i}{x} \approx \pdv[2]{\Tilde{u}_i}{x} = \sum_{j \in \mathcal{C}_i} \lambda_{ij} (u_j - u_i). \label{eq:temp}
\end{align} Using the first order derivative in \eqref{eq:temp}, we obtain a semi-discretized form of the the linear advection equation \eqref{eq:linearAdvection}. Rather than directly use the neighbouring points \(u_j\), we again use a midpoint value \(u_{ij}\) that we define below. To compensate for the smaller stencil and allow for the reuse of the coefficients \(\alpha_{ij}\), we must add a factor two. This yields the following semi-discretized scheme \begin{align}
	\pdv{u_i}{t} &= - 2a \sum_{j \in \mathcal{C}_i} \kappa_{ij} (u_{ij} - u_i). \label{eq:1DMUSCLGradient}
\end{align} The selection of the midpoint is done in an upwind manner \begin{align}
	u(\frac{x_i + x_j}{2}) \approx u_{ij} &= \begin{cases}
		\Bar{u}_{ij}, & \text{if } a \Delta x_{ij} > 0 \\
		\Bar{u}_{ji}, & \text{else }
	\end{cases},
\end{align} where $\Bar{u}_{ij}$ and $\Bar{u}_{ji}$ are quadratic reconstructions using a Taylor expansion from $u_i$ and $u_j$ to the midpoint \begin{align}
	u(\frac{x_i + x_j}{2}) &= u_i + \frac{\Delta x_{ij}}{2}\pdv{u_i}{x} + \frac{\Delta x_{ij}^2}{8}\pdv[2]{u_i}{x} + \mathcal{O}(\Delta x_{ij}^3) \\
	\approx \Bar{u}_{ij} &= u_i + \frac{\Delta x_{ij}}{2} \sum_{k \in C_i} \kappa_{ik}(u_k - u_i) + \frac{\Delta x_{ij}^2}{8} \sum_{k \in C_i} \lambda_{ik}(u_k - u_i), \label{eq:1DMUSCLmidpointRec1} \\
	u(\frac{x_i + x_j}{2}) &= u_j - \frac{\Delta x_{ij}}{2}\pdv{u_j}{x} + \frac{\Delta x_{ij}^2}{8}\pdv[2]{u_j}{x} + \mathcal{O}(\Delta x_{ij}^3) \\
	\approx \Bar{u}_{ji} &= u_j - \frac{\Delta x_{ij}}{2} \sum_{k \in C_j} \kappa_{jk}(u_k - u_j) + \frac{\Delta x_{ij}^2}{8} \sum_{k \in C_j} \lambda_{jk}(u_k - u_j). \label{eq:1DMUSCLmidpointRec2}
\end{align} The scheme above is of second order in space, or one order lower than the order of the Taylor expansion that is computed for the least square problem. Extension to higher order is straightforward: solve the least square problem for arbitrary order $m$, use all derivatives available for the reconstruction at all midpoints, and finally, compute the first order derivative using the midpoints as in \eqref{eq:1DMUSCLGradient}. As for the time integration method, we chose explicit Runge-Kutta methods. 

By selecting the stencil with which to reconstruct the midpoint in an upwind manner, the final scheme has a stencil weighted in the upwind direction and is, therefore, $L^2$ stable. This is verified numerically in section \ref{section:Stability}. In addition, the scheme is efficient because it requires only one central stencil that is reused in the computation of the reconstruction.

Although the scheme above is stable in the \(L^2\) norm, it is possible to show that it cannot be positive (stable in the \(L^{\infty}\) norm). To avoid spurious oscillations one could introduce some kind of limiting as in \cite{chandrashekar_positive_2004, huh_new_2018}. Instead, we opt for the more simple MOOD method. In the MOOD method, so-called `admissibility criteria' are explicitly checked at each time step. If the numerical solution does not meet the criteria, the order of the spatial discretisation method is reduced. Admissibility criteria typically include the discrete maximum principle (DMP) or a relaxed version of it \cite{diot_improved_2012} \begin{align}
	\min_{j \in C_i}(u_i^n, u_j^n) \leq u_i^{n+1} \leq \max_{j \in C_i}(u_i^n, u_j^n). \label{eq:DMP}
\end{align} By explicitly enforcing the DMP, positivity of the solution is ensured and thus spurious oscillations at discontinuities are avoided. There exist several strategies with which to decrease the order of the method. In this text, upon failure of the admissibility criterion, the spatial order is immediately reduced to order one. For the first order scheme, we use the upwind scheme from section \ref{section:FirstOrder}, which is positive under a known CFL condition. This MOOD procedure is then embedded into the stages of an explicit Runge-Kutta time integration scheme. At a MOOD event, the Runge-Kutta substage is replaced from a forward Euler step to the next Runge-Kutta stage. The resulting behaviour of MOOD method is that the scheme drops to order one around discontinuities and achieves high order in smooth regions. 

A strict discrete maximum property such as \eqref{eq:DMP} was also used in the original work on MOOD \cite{clain_high-order_2011}. Indeed this makes sense for the linear advection equation, since all scalar hyperbolic equations satisfy the discrete maximum property \cite{bressan_hyperbolic_2005}. However, in \cite{clain_high-order_2011, diot_improved_2012} it was found, that if the discrete maximum property is enforced with MOOD, the order of the scheme is limited to two. To circumvent this issue, multiple relaxed formulations of the DMP criterion have been proposed. In \cite{farmakis_weno_2020, fernandez-fidalgo_posteriori_2018, loubere_cat-mood_2024}, the left and right hand side of the DMP inequality \eqref{eq:DMP}, are relaxed with a parameter $\delta$. The parameter, if chosen large enough, can recover the order of the higher-order scheme, but can also potentially reintroduce oscillations at discontinuities. In addition, this parameter inevitably depends on the grid and problem and therefore requires some tuning. In \cite{diot_multidimensional_2013, diot_improved_2012}, it was observed that the loss of order of the original MOOD method occurs at local extrema. This prompted us to relax the DMP criterion at local extrema. Local extrema are detected if the signs of the minimal and maximal curvature are equal. Local extrema are then considered actual extrema, and not spurious oscillations or discontinuities, based on a smoothness criterion. This relaxed MOOD criterion is referred to as the \emph{u2 detection procedure} and is also used in \cite{boscheri_direct_2015, bourgeois_gp-mood_2022}. To our understanding, the smoothness criterion on the curvatures in Definition 3.4 in \cite{diot_multidimensional_2013} is ill-defined for local maxima (negative curvatures) because the ratio of curvatures does not lie in $[0, 1]$. The u2 detection procedure used in this text corrects for this, and is outlined below. Define \begin{align}
	\Tilde{\mathcal{X}}_i^{min} &= \min_{j \in C_i} \left( \left| \pdv[2]{\Tilde{u}_i}{x} \right| , \left| \pdv[2]{\Tilde{u}_j}{x} \right| \right) &&\text{and} \; \Tilde{\mathcal{X}}_i^{max} = \max_{j \in C_i} \left( \left| \pdv[2]{\Tilde{u}_i}{x} \right|, \left| \pdv[2]{\Tilde{u}_j}{x} \right| \right) \\
	\mathcal{X}_i^{min} &= \min_{j \in C_i} \left( \pdv[2]{\Tilde{u}_i}{x}, \pdv[2]{\Tilde{u}_j}{x} \right)  &&\text{and} \; \mathcal{X}_i^{max} = \max_{j \in C_i} \left( \pdv[2]{\Tilde{u}_i}{x}, \pdv[2]{\Tilde{u}_j}{x} \right)
\end{align} The improved u2 detection criterion then consists of checking the DMP \eqref{eq:DMP}. A solution that does not satisfy 
the DMP, is still eligible if it satisfies \begin{align}
	\mathcal{X}_i^{min} \mathcal{X}_i^{max} > 0, \label{eq:MOODu3A}\\
	\frac{\Tilde{\mathcal{X}}_i^{min}}{\Tilde{\mathcal{X}}_i^{max}} \geq 1/2 \label{eq:MOODu3B}. 
\end{align}

Condition \eqref{eq:MOODu3A} determines whether the numerical solution can be considered non-oscillating. Condition \eqref{eq:MOODu3B} determines if the numerical solution is sufficiently smooth. Finally, as was observed in \cite{diot_multidimensional_2013}, small micro-oscillations in a flat region of the solution can falsely activate the curvature criteria. To make the detection criteria robust to micro-oscillations, all criteria are relaxed with a small $\delta$. Thus, a solution that does not adhere to the DMP criterion \eqref{eq:DMP} is still accepted if \begin{align} 
	\left| \max_{j \in C_i}(u_i^n, u_j^n) - \min_{j \in C_i}(u_i^n, u_j^n) \right| \leq \delta^3. \label{eq:RelaxedDMP}
\end{align} If a solution does not satisfy condition \eqref{eq:RelaxedDMP}, then the relaxed u2 detection procedure is checked. \begin{align}
	\mathcal{X}_i^{min} \mathcal{X}_i^{max} > - \delta, \label{eq:MOODu4A} \\
	\left( \frac{\Tilde{\mathcal{X}}_i^{min}}{\Tilde{\mathcal{X}}_i^{max}} \geq 1/2 \right) \ \text{or} \ \left( \Tilde{\mathcal{X}}_i^{max} < \delta \right) \label{eq:MOODu4B}. 
\end{align} The additional check in \eqref{eq:MOODu4B} compared to \eqref{eq:MOODu3B} is to ensure that the curvature ratio criterion is only checked if the curvatures are significantly large, i.e., the solution is not flat. The parameter $\delta$ is defined as the size of the cell, as in \cite{diot_multidimensional_2013}.

\subsubsection{Formulation in two space dimensions}
\label{section:MUSCLMOOD2D}
Consider a third order MLS approximation in two spatial dimensions. This procedure yields a set of coefficients \(\alpha_{ij}, \beta_{ij}, \gamma_{ij}, \eta_{ij}\) and \(\nu_{ij}\) with which the following derivatives can be approximated
\begin{align}
	\pdv{u_i}{x} &\approx \pdv{\Tilde{u}_i}{x} = \sum_{j \in \mathcal{C}_i} \alpha_{ij} (u_j - u_i), \; \pdv{u_i}{y} \approx \pdv{\Tilde{u}_i}{y} = \sum_{j \in \mathcal{C}_i} \beta_{ij} (u_j - u_i), \label{eq:approx1} \\ 
	\pdv{u_i}{x}{y} &\approx \pdv{\Tilde{u}_i}{x}{y} = \sum_{j \in \mathcal{C}_i} \gamma_{ij} (u_j - u_i), \; \pdv[2]{u_i}{x} \approx \pdv[2]{\Tilde{u}_i}{x} = \sum_{j \in \mathcal{C}_i} \eta_{ij} (u_j - u_i), \\
	\pdv[2]{u_i}{y} &\approx \pdv[2]{\Tilde{u}_i}{y} = \sum_{j \in \mathcal{C}_i} \nu_{ij} (u_j - u_i).
\end{align} As before, we discretize the divergence using approximations of the derivatives obtained with MLS \eqref{eq:approx1} \begin{align}
	\pdv{u_i}{t} &= - 2 a_x \sum_{j \in C_i} \alpha_{ij} (u_{ij} - u_i) - 2a_y \sum_{j \in C_i} \beta_{ij} (u_{ij} - u_i),
\end{align} and define the midpoint $u_{ij}$ in an upwind manner\begin{align}
	u_{ij} &= \begin{cases}
		\Bar{u}_{ij}, & \text{if } a_x \Delta x_{ij} + a_y \Delta y_{ij} > 0 \\
		\Bar{u}_{ji}, & \text{else. } 
	\end{cases}
\end{align} The reconstructions from the centre point and the neighbouring point are computed using a Taylor series in which we have substituted the exact derivatives for their MLS approximations \begin{align}
	\Bar{u}_{ij} &= u_i + \frac{\Delta x_{ij}}{2}\pdv{\Tilde{u}_i}{x} + \frac{\Delta y_{ij}}{2}\pdv{\Tilde{u}_i}{y} + \frac{\Delta x_{ij}^2}{8}\pdv[2]{\Tilde{u}_i}{x} + \frac{\Delta y_{ij}^2}{8}\pdv[2]{\Tilde{u}_i}{y} + \frac{\Delta x_{ij} \Delta y_{ij}}{4} \pdv{\Tilde{u}_i}{x}{y}, \label{eq:2DMUSCLmidpointRec1} \\
	\Bar{u}_{ij} &= u_j - \frac{\Delta x_{ij}}{2}\pdv{\Tilde{u}_j}{x} - \frac{\Delta y_{ij}}{2}\pdv{\Tilde{u}_j}{y} + \frac{\Delta x_{ij}^2}{8}\pdv[2]{\Tilde{u}_j}{x} + \frac{\Delta y_{ij}^2}{8}\pdv[2]{\Tilde{u}_j}{y} + \frac{\Delta x_{ij} \Delta y_{ij}}{4} \pdv{\Tilde{u}_j}{x}{y} . \label{eq:2DMUSCLmidpointRec2}
\end{align} As before, this method can be coupled with MOOD to obtain a positive scheme. The additional conditions for the curvatures in the u2 detection procedure are checked in both spatial directions.

\subsection{Meshless WENO method}
\label{section:WENO}
\subsubsection{Formulation in one space dimension}
\label{section:WENO1D}
Due to the efficiency, accuracy and sharp behaviour at discontinuities, Weighted Essentially Non-Oscillatory methods are very popular. In \cite{avesani_new_2014}, a moving least square method was combined with a WENO scheme in the context of SPH-type particle methods. In \cite{tiwari_meshfree_2022}, this MLS-WENO method was applied to a meshless method for the BGK-Boltzmann equation. We summarise the method from \cite{tiwari_meshfree_2022} below and use it as a benchmark for the newly developed MUSCL-like schemes. We again write out the second order scheme and discuss generalisations to higher orders. 

Let $L_i$ and $R_i$ be the sets of points containing the points to the left and right of grid point $i$, such that $C_i = L_i \cup R_i$. Considering a third-order Taylor expansion and by applying the least squares method from \ref{section:leastSquares}, we obtain first order and second order derivatives for each neighbourhood, e.g. $\pdv{u_{i, L_i}}{x}$. The first order derivative that is ultimately used in the time integration routine is the weighted sum \begin{align}
	\pdv{u_i}{x} = \omega_{L_i} \pdv{u_{i, L_i}}{x} + \omega_{C_i} \pdv{u_{i, C_i}}{x} + \omega_{R_i} \pdv{u_{i, R_i}}{x}.
\end{align} The weights $\omega_k$ are defined by \begin{align}
	\omega_k = \frac{\beta_k}{\beta_{C_i} + \beta_{R_i} + \beta_{L_i}}, \; \; k = R_i, L_i, C_i,
\end{align} with \begin{align}
	\beta_k = \frac{D_k}{\left( \pdv{u_{i, k}}{x}^2 \Delta x^2 + \pdv[2]{u_{i, k}}{x}^2 \Delta x^4 + \varepsilon \right)^2 }, \; \; k = R_i, L_i, C_i. \label{eq:1DWENOCoeff}
\end{align} The parameter $\varepsilon$ is set to a small value to avoid division by zero, and $\Delta x$ is the initial average grid point spacing. The coefficients $D_k$ are chosen in an upwind manner \begin{align}
	\begin{cases}
		D_{L_i} = 0.5, D_{C_i} = 0.5, D_{R_i} = 0 \; \; \text{if } a \geq 0 \\
		D_{L_i} = 0, D_{C_i} = 0.5, D_{R_i} = 0.5 \; \; \text{else. }
	\end{cases}
\end{align}
The WENO method selects the most well-behaved stencil based on the size of the derivatives. In a smooth region, the central stencil will most often be selected, around a discontinuity, the relevant one-sided stencil is used. The upwind procedure avoids obvious stability issues. In theory, this method is easily extended to higher orders: one can simply use higher order Taylor expansions and include the higher order derivatives in the denominator of \eqref{eq:1DWENOCoeff}. Note that compared to the MUSCL method from section \ref{section:MUSCLMOOD1D}, per point two least squares problems must be solved, instead of one.

\subsubsection{Formulation in two space dimensions}
\label{section:WENO2D}
The extension of the WENO method to two dimensions is straightforward. We divide the neighbours of a particle into subsets left $L_i$, right $R_i$, top $T_i$ and bottom $B_i$, with the central stencil $C_i = L_i \cup R_i \cup T_i \cup B_i$. For all five stencils, the first and second order derivatives can be computed using the procedure from section \ref{section:leastSquares}. The unnormalised WENO coefficients $\beta_k$ are then computed by \begin{align}
	\beta_k = \frac{D_k}{\left( \pdv{u_{i, k}}{x}^2 \Delta x^2 + \pdv[2]{u_{i, k}}{x}^2 \Delta x^4 + \pdv{u_{i, k}}{y}^2 \Delta y^2 + \pdv[2]{u_{i, k}}{y}^2 \Delta y^4 + \pdv{u_{i, k}}{x}{y}^2 \Delta x^2 \Delta y^2 + \varepsilon \right)^2}, 
\end{align} with $k$ either $L_i, C_i$ or $R_i$ for the derivative in the x-direction, or $D_i, C_i$ or $U_i$ for the derivative in the y-direction. The coefficient $D_k$ is chosen such that the upwind and central stencil has a weight of 0.5, and the `downwind' stencil has a weight of 0. The weights $\omega_k$ and the final derivative used by the WENO scheme are computed analogously as the one-dimensional case. 

\begin{remark}
	The WENO scheme used here is based on the scheme from \cite{avesani_new_2014}. The scheme from \cite{avesani_new_2014} proved to be very unstable due to the large weight applied to the central stencil and non-zero weights for the `downwind' stencils.
\end{remark}

\section{Numerical Experiments}
\label{section:Numerics}
This section aims to assess the performance of the newly developed MUSCL-like scheme. To this end, all schemes were implemented in a Julia package \textit{meshfree4hypeq} \cite{githubRepo}. In section \ref{section:Convergence}, we do a convergence analysis. In section \ref{section:Stability}, we asses the linear stability of the MUSCL scheme. Due to the irregularity of the grid, a classical Von Neumann stability analysis is not possible and we must resort to numerical calculations of spectra. As previously remarked, the schemes used here are not conservative. To verify that the methods remain robust when Dirichlet conditions are used, we examine this case in section \ref{section:DBoundary}. In section \ref{section:Conservation}, we illustrate the lack of conservation of the MUSCL method compared to a classical upwind meshless method and the WENO method. Finally in section \ref{section:Efficiency}, we compare the efficiency of the MUSCL scheme and the WENO scheme in obtaining an error for a predefined computation time. 

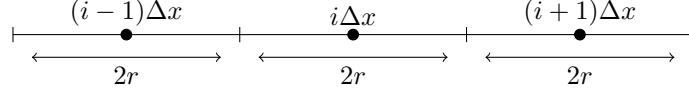
\begin{figure}
	\centering
	\begin{tikzpicture}
		\draw[-] (-4.5,0) -- (4.5,0) node[right] {};
		
		\filldraw ( -3, 0) circle (2pt) node[above] {$(i-1)\Delta x$};
		\filldraw (  0, 0) circle (2pt) node[above] {$i\Delta x$};
		\filldraw (  3, 0) circle (2pt) node[above] {$(i+1)\Delta x$};
		
		\draw[<->] (-1.25, -0.3) -- (1.25, -0.3) node[midway, below] {$2r$};
		\draw[<->] (-4.25, -0.3) -- (-1.75, -0.3) node[midway, below] {$2r$};
		\draw[<->] (4.25, -0.3) -- (1.75, -0.3) node[midway, below] {$2r$};
		
		\draw (-1.5, 0.1) -- (-1.5, -0.1);
		\draw (1.5, 0.1) -- (1.5, -0.1);
		\draw (4.5, 0.1) -- (4.5, -0.1);
		\draw (-4.5, 0.1) -- (-4.5, -0.1);
	\end{tikzpicture}
	
	\caption{Illustration of how random grids are generated with $r \in \left[0, \frac{\Delta x}{2} \right]$.}
	\label{fig:1DGridGeneration}
\end{figure}

In all the numerical experiments that follow, we consider square domains. The irregular grids are generated from a uniform grid with spacing $\Delta x$ with additional noise as \begin{align}
	x_i = \Delta x(i-1) + r(2p-1), \; \; \text{with } p \sim U(0, 1). \label{eq:randomness}
\end{align} The parameter $r$ is referred to as the `randomness' and set to $\frac{\Delta x}{2}$ unless stated otherwise. Also, see figure \ref{fig:1DGridGeneration}. To generate a grid in 2D, the same methodology is applied as in 1D. Other globally fixed parameters for the 1D and 2D simulations are summarized in table \ref{tab:1DGlobalParams}.

\begin{table}[h!]
	\centering
	\begin{tabular}{|c||c|c|c|c|c|} 
		\hline
		& $a$ & $h_{max}$ & $\varepsilon$ & \(\alpha \)\\ [0.5ex] \hline \hline
		1D & 1 & $3.5 \Delta x$ & $1\times10^{-6}$ & \( \Delta x^{-2} \) \\
		2D & (1, 1) & $\sqrt{34} \Delta x$ & $1\times10^{-12}$ & \( \frac{6}{h_{max}^2} \) \\
		\hline
	\end{tabular}
	\caption{Global fixed parameters used in the simulations. From left to right the parameters are the velocity in the linear advection equation $a$, the maximum neighbour distance $h_{max}$, the WENO parameter $\varepsilon$, and the dimension-dependent weight function parameter \(\alpha\) (see \eqref{eq:WeightFunction}. }
	\label{tab:1DGlobalParams}
\end{table}

\subsection{Convergence}
\label{section:Convergence}

\subsubsection{One space dimension}
\begin{figure}[h!]
	\centering
	\begin{subfigure}[b]{0.49\textwidth}
		\centering
		\includegraphics[width=\linewidth]{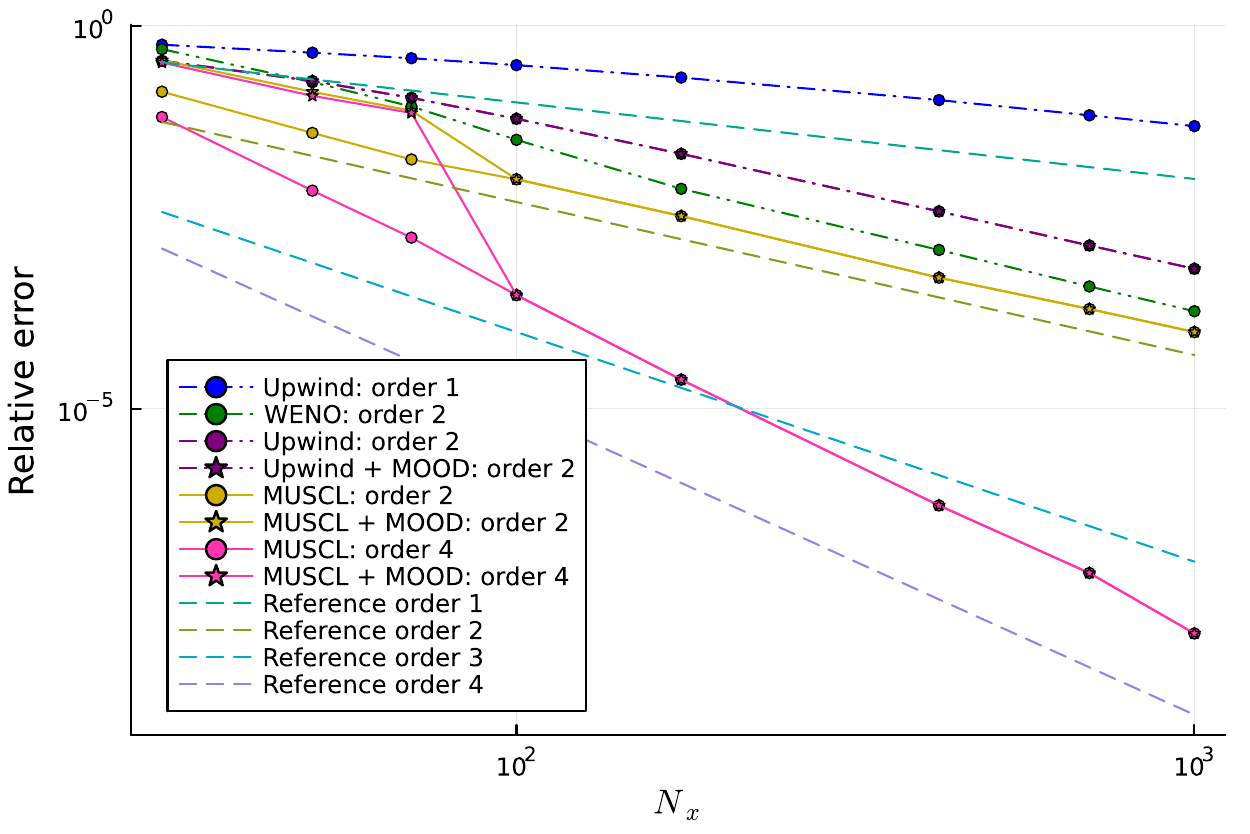}
		\label{fig:1DConvergenceSmooth}
	\end{subfigure}
	\begin{subfigure}[b]{0.49\textwidth}
		\centering
		\includegraphics[width=\linewidth]{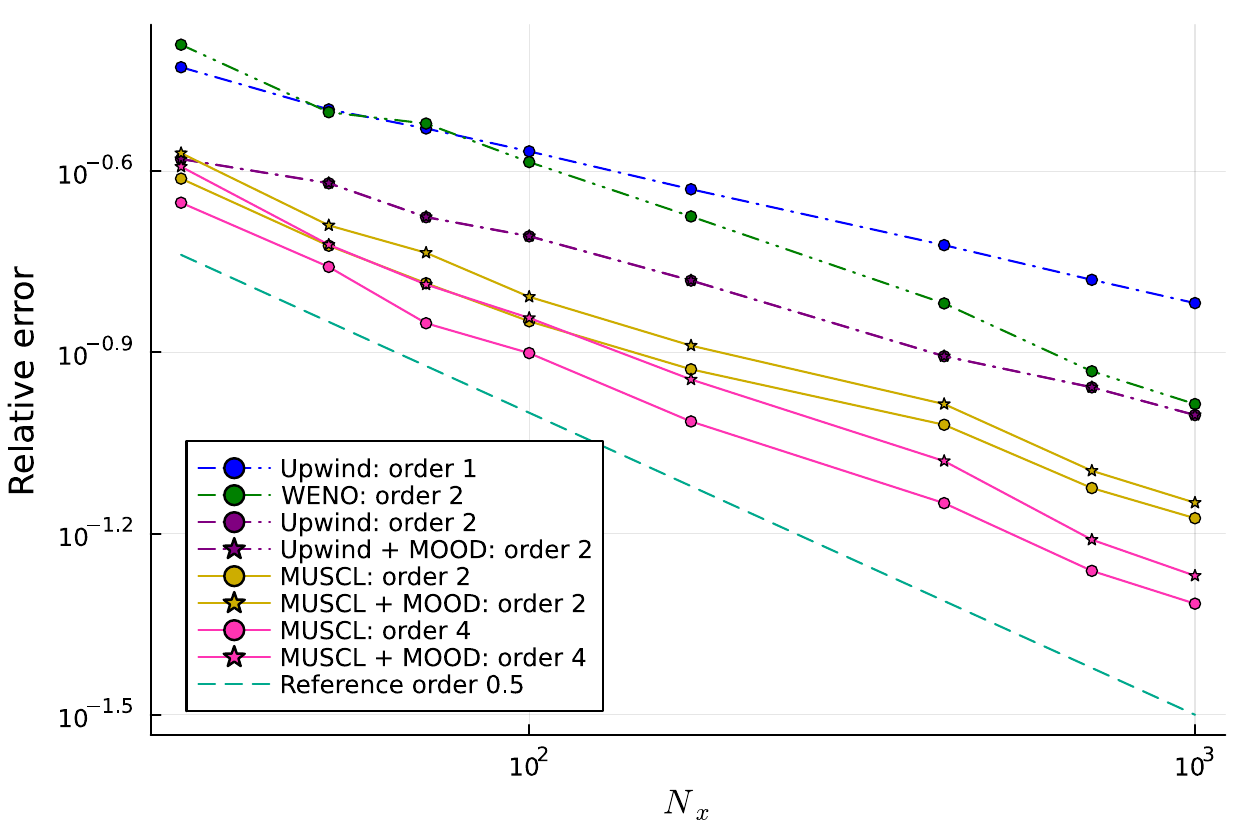}
		\label{fig:1DConvergenceShock}
	\end{subfigure}
	\caption{Convergence of 1D meshless schemes for smooth initial condition (left) and shock initial condition (right). The $x$-axis represents the amount of grid points. The $y$-axis plots the relative error. The numbers in the legend are the theoretical orders of the methods.}
	\label{fig:1DConvergence}
\end{figure}

In this section, we show the convergence of the MUSCL and WENO schemes described in section \ref{section:HigherOrder} in one spatial dimension. We consider a periodic domain $\Omega = \left[-5, 5\right]$ on which we solve the linear advection equation with initial conditions \begin{align}
	u_1(0, x) = \exp \left( -x^2 \right) \; \text{and} \; \; u_2(0, x) = \begin{cases}
		1 \; \; \text{if } x > 0 \\
		0 \; \; \text{else.}
	\end{cases} \label{eq:initialConditions1D}
\end{align} up to time $t=2.5$. The numerical solution is computed with two upwind schemes, a WENO scheme and two MUSCL schemes, all combined with a third-order Runge-Kutta time integrator. The time step used for all methods is the time step for which the first order scheme is positive \eqref{eq:1DCFLCondition}, with a small CFL number $1/20$ such that the time integration error is not visible. The error is calculated by computing the 2-norm of the error: the exact solution on the irregular grid and the numerical solution. Finally, the results are plotted in figure \ref{fig:1DConvergence}. 

\begin{figure}[h!]
	\centering
	\begin{subfigure}[b]{0.49\textwidth}
		\centering
		\includegraphics[width=\linewidth]{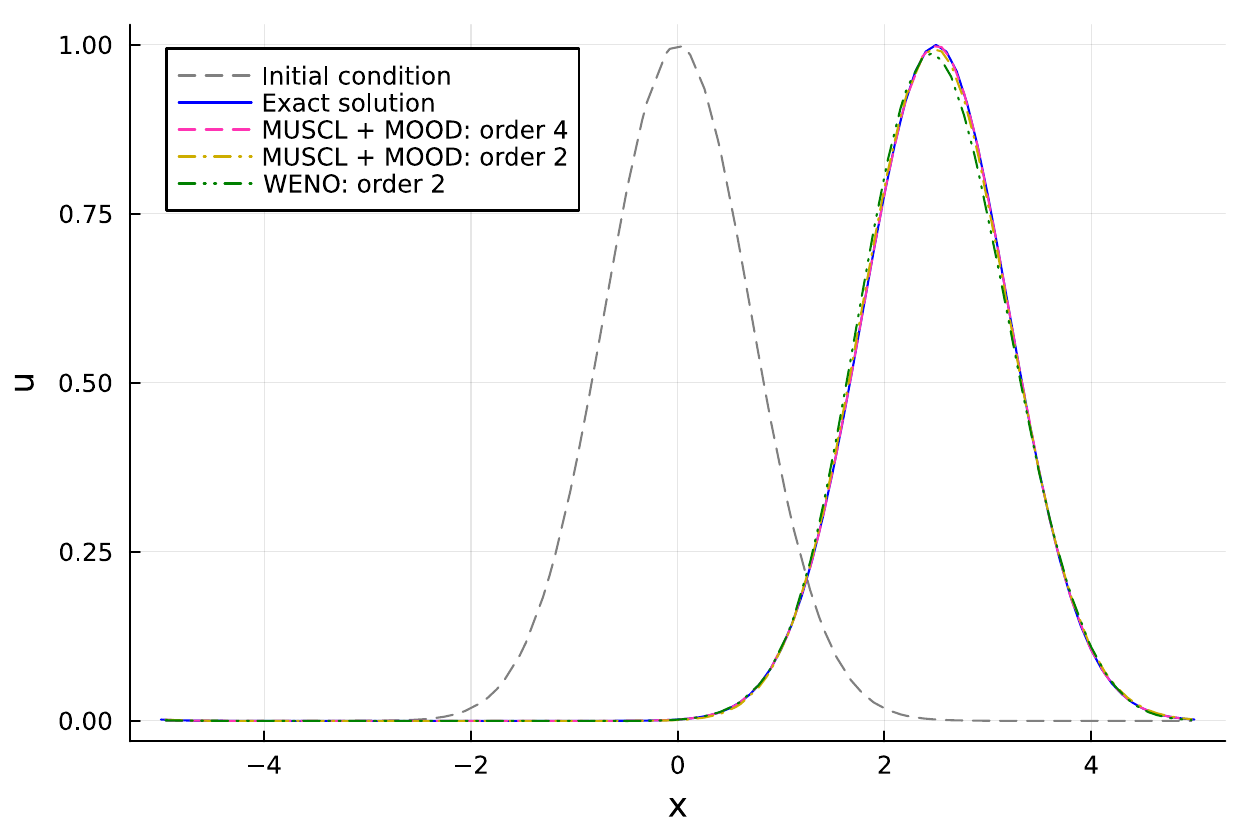}
		\label{fig:1ConvergenceSmoothSolution}
	\end{subfigure}
	\begin{subfigure}[b]{0.49\textwidth}
		\centering
		\includegraphics[width=\linewidth]{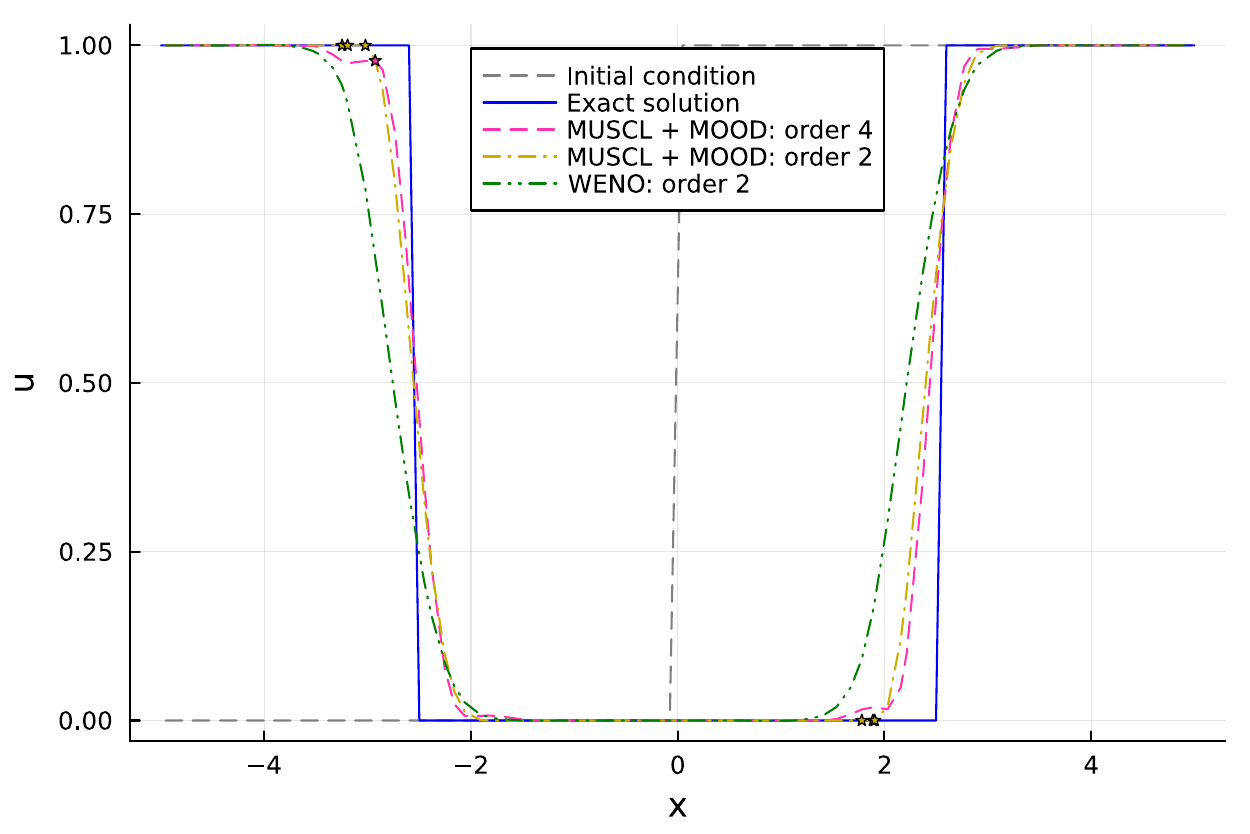}
		\label{fig:1ConvergenceShockSolution}
	\end{subfigure}
	\caption{Solution for the sine and shock initial condition for several schemes with $N_x = 100$. Stars indicate the location of a MOOD event at the final time step.}
	\label{fig:1ConvergenceSolutionPlots}
\end{figure}

 All 1D schemes achieve the order as advertised. We achieved up to fourth order using the MUSCL reconstruction technique. From the figure, it follows that the improved u2 detection procedure maintains the order of the schemes, granted that the grids are sufficiently fine. The first order and third-order MUSCL schemes are missing from the convergence plot because they are unstable (see section \ref{section:Stability}). In figure \ref{fig:1ConvergenceSolutionPlots}, the solutions are plotted for 100 grid points. Both from the convergence plot and from the solution plots it can be observed that the second order MUSCL scheme slightly outperforms the WENO scheme. The fourth order MUSCL scheme yields the best results. Due to the relaxation of the MOOD criterion at local extrema, the MOOD schemes capture the peak of the Gaussian well. In case of the shock solution, the MOOD scheme yields minimal MOOD events around the discontinuity. In addition, no MOOD events are observed in the region where the solution is flat, due to the relaxation of the MOOD parameter $\delta$.

\subsubsection{Two space dimensions}

\begin{figure}[h!]
	\centering
	\begin{subfigure}[b]{0.49\textwidth}
		\centering
		\includegraphics[width=\linewidth]{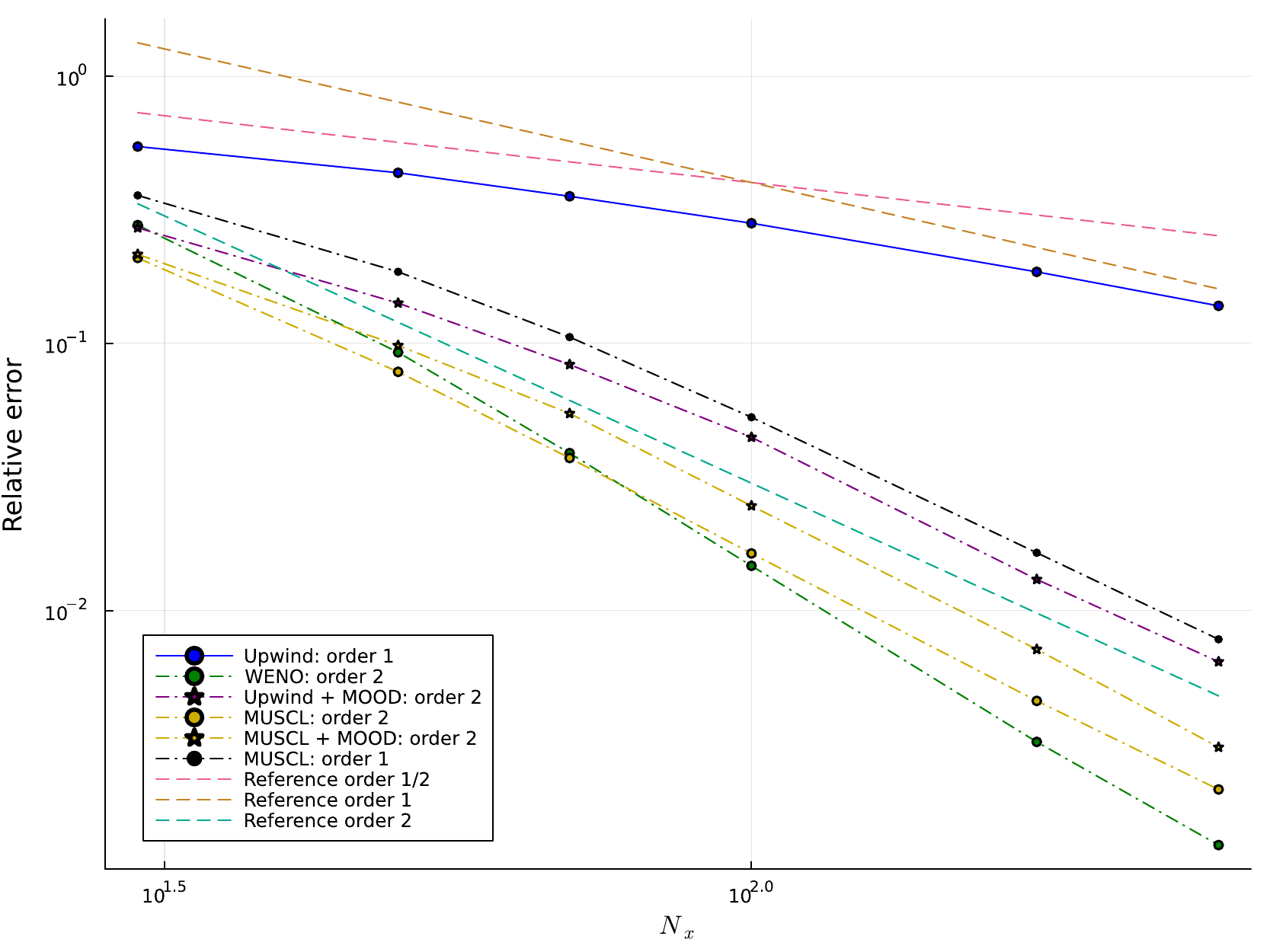}
		\label{fig:2DConvergenceSmooth}
	\end{subfigure}
	\begin{subfigure}[b]{0.49\textwidth}
		\centering
		\includegraphics[width=\linewidth]{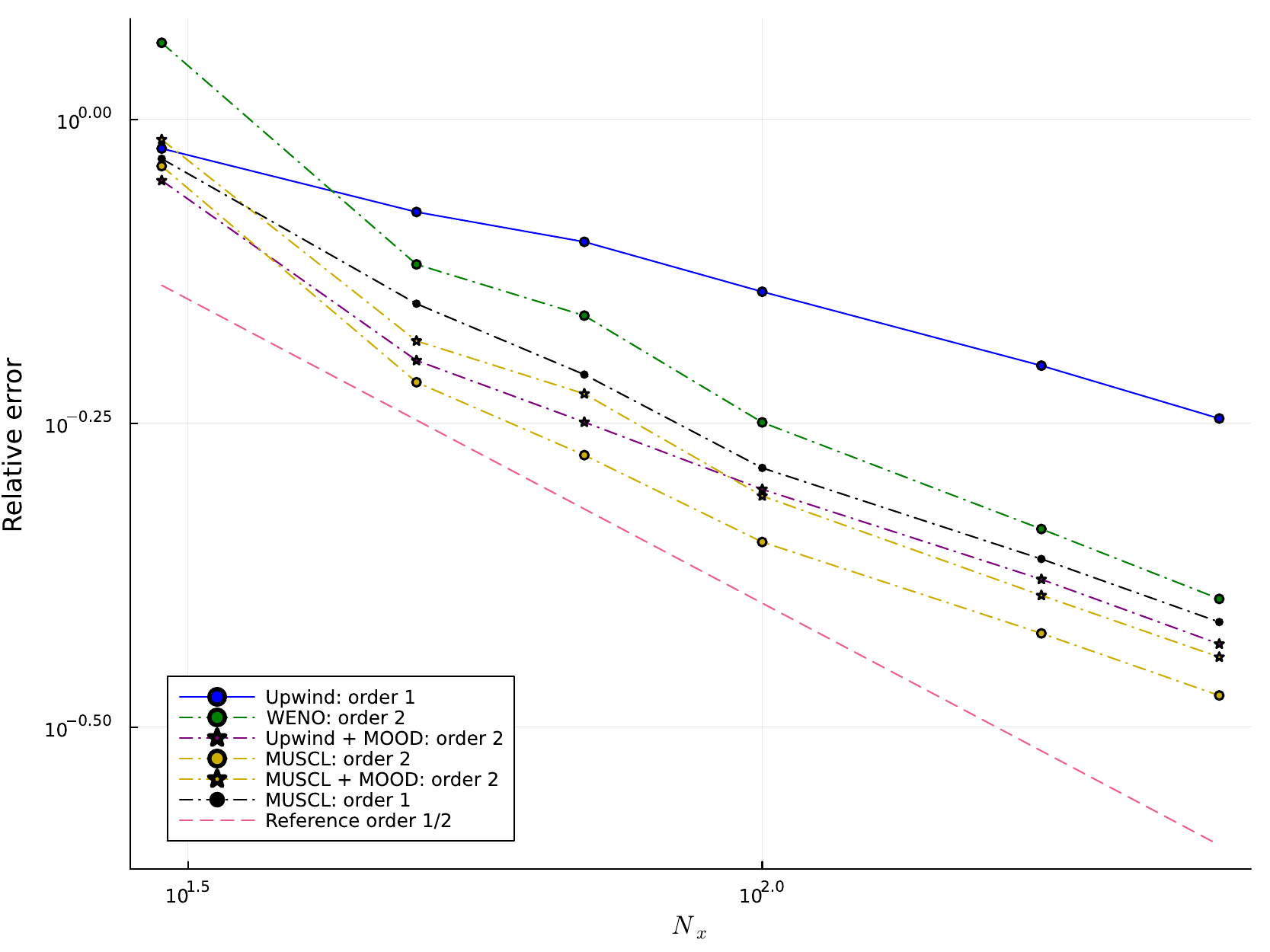}
		\label{fig:2DConvergenceShock}
	\end{subfigure}
	\caption{Convergence of 2D meshless schemes for smooth initial condition (left) and shock initial condition (right). The $x$-axis represents the amount of grid points. The $y$-axis plots the relative error. The numbers in the legend are the theoretical orders of the methods.}
	\label{fig:2DConvergence}
\end{figure}

We repeat the analysis in two spatial dimensions. We consider a 2D periodic domain \( [-5, 5]^2 \) on which we solve the linear advection equation with initial conditions \begin{align}
	u_1(0, x, y) = \exp(-x^2 - y^2) \; \text{and} \; u_2(0, x, y) = \begin{cases*}
		1 \; \; \text{if } -0.5 < x, y < 0.5 \\
		0 \; \; \text{else.}
	\end{cases*} \label{eq:2DSmoothShockInit}
\end{align} up to time t = \(1.0\). The numerical solution is computed using the second order WENO scheme, the second order MUSCL scheme, the first and second order MUSCL scheme with MOOD, the second order upwind scheme, and the first order scheme from section \ref{section:FirstOrder2D}. All second order schemes are paired with a third-order Runge-Kutta time integration scheme. The time step is chosen such that the first order scheme is stable with CFL = \(1/40\). The error is computed as the 2-norm of the numerical solution minus the exact solution on the irregular grid. The results are plotted in figure \ref{fig:2DConvergence}. All schemes achieve their advertised order. Surprisingly, the first order MUSCL scheme obtains second order convergence. The upwind scheme from section \ref{section:FirstOrder2D} is in theory a first order method. In practice, the scheme only achieves this order for very fine grids. This is not entirely unusual due to the additional diffusion that is added to the scheme to keep it positive. 

\subsection{Stability \& sensitivity}
\label{section:Stability}

\subsubsection{One space dimension}
In this section, we analyze the stability of the meshless schemes by linear matrix stability. Due to the non-linear nature of the WENO schemes, their stability can only be analysed through direction simulation. These schemes are therefore not analyzed here. Instead, we consider upwind schemes, central schemes and MUSCL schemes of various orders without the MOOD procedure. The central scheme directly uses the derivative obtained from the least squares procedure using a central stencil. We proceed as follows. First, we consider one irregular grid $\Omega = \left[ -5, 5\right]$ with 100 grid points and plot the spectra of the semi-discretised ODEs. We then repeat this process for many different grids to assess how sensitive the stability of the schemes is to different grids. 

\begin{figure}[h!]
	\centering
	\begin{subfigure}[b]{0.49\textwidth}
		\centering
		\includegraphics[width=\linewidth]{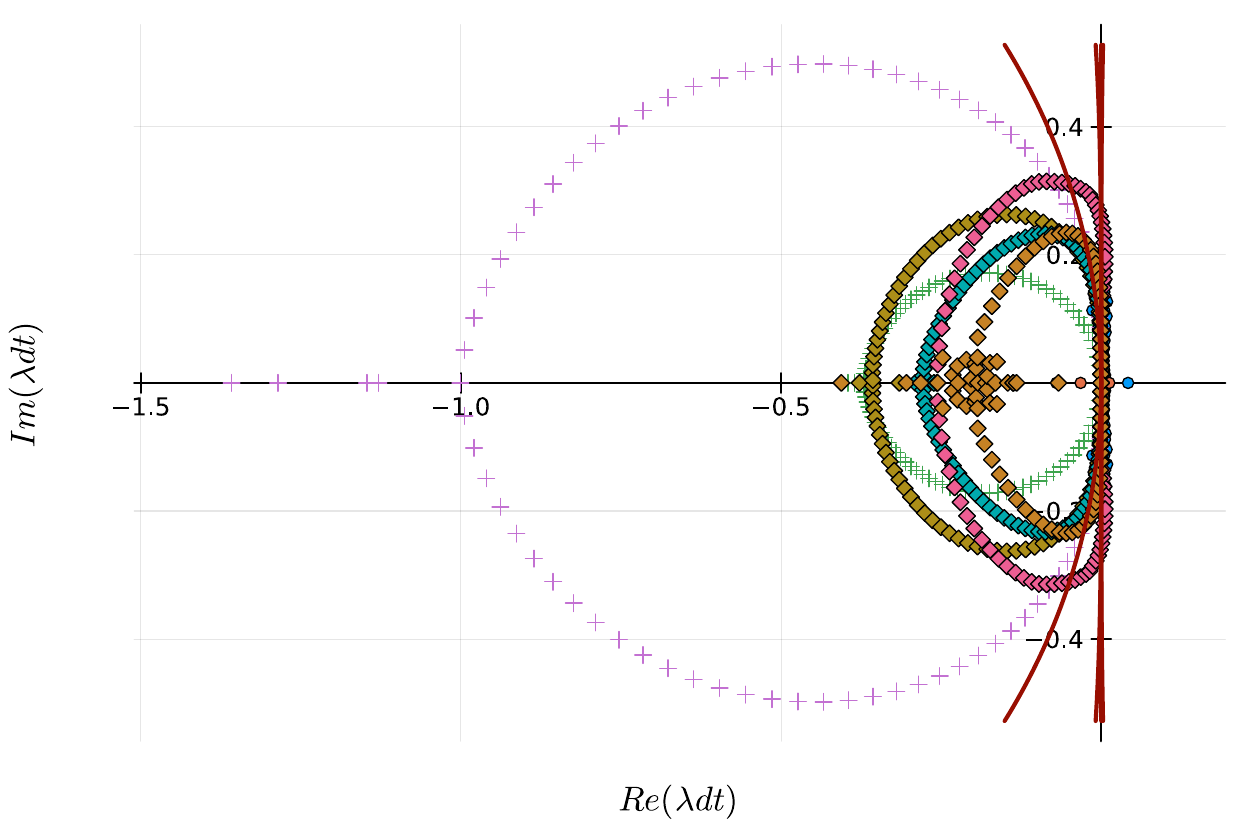}
		\label{fig:1DEigenvaluesFull}
	\end{subfigure}
	\begin{subfigure}[b]{0.49\textwidth}
		\centering
		\includegraphics[width=\linewidth]{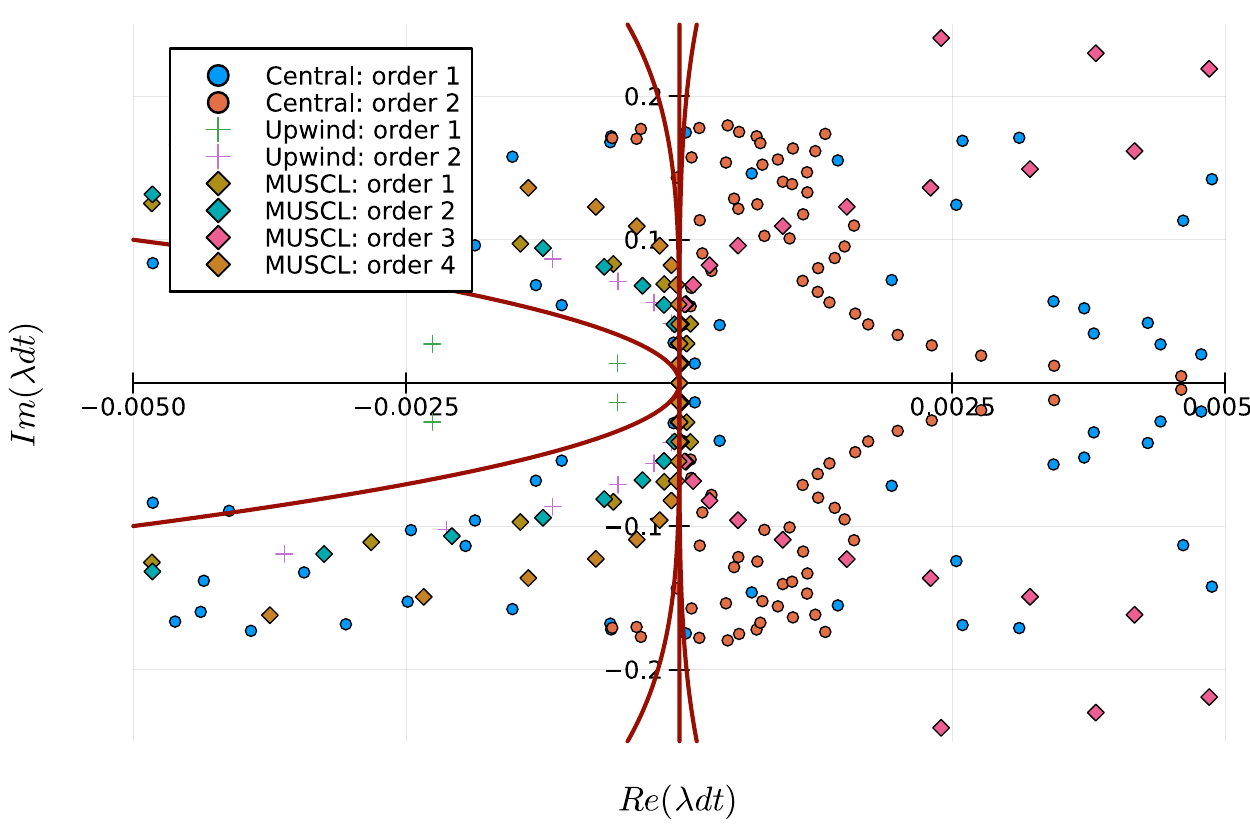}
		\label{fig:1DEigenvaluesZoom}
	\end{subfigure}
	\caption{Spectra of several meshless schemes with CFL number $1/8$. The lines in the graph are the stability regions for the first, second, third and fourth order explicit Runge-Kutta method.}
	\label{fig:1DEigenvalues}
\end{figure}

In figure \ref{fig:1DEigenvalues}, the spectra of several meshless schemes are plotted. As expected the upwind schemes yield stable spectra. In the case of the central scheme, the eigenvalues fall to the left and right of the imaginary axis, yielding an unstable scheme. In the fixed grid case, these eigenvalues would lie exactly on the imaginary axis, such that a third or fourth order Runge-Kutta scheme can stabilize the scheme. This is however not the case for irregular grids, and as a result, standard GFDMs cannot be directly applied to hyperbolic equations. The first order and third-order MUSCL schemes consistently have multiple eigenvalues in the right half plane. The second order and fourth order schemes have spectra that lie fully in the left-half plane. As the order of the MUSCL schemes increases, the spectra tend to be smaller and `stick' much longer to the imaginary axis. Higher-order Runge-Kutta methods that include the imaginary axis in their stability domain are thus still relevant for the MUSCL schemes. We note that in the case of a uniform grid, all MUSCL schemes, as well as the third-order and first order schemes, become stable. Stability properties can then be proven using a classical von Neumann analysis. 

\begin{remark}
	The unstable schemes can be stabilized by adding a small viscosity term that shifts the eigenvalues into the left half-plane.
\end{remark}

It is not sufficient to analyze the stability of the schemes for just one irregular grid alone, since already for the first order scheme the stability condition depends on the position of every point \eqref{eq:1DCFLCondition}. We therefore repeat the analysis above for 100 different randomly generated irregular grids with varying amount of grid points. We then log if a scheme is unstable for a particular grid if there is at least one eigenvalue with a real part larger than $1e-13$. Since these schemes are intended to be used in moving mesh methods, we can also take into account that there is some minimum distance between grid points that is always satisfied. This is typically implemented in a grid management routine, see \cite{tiwari_meshfree_2022, suchde_conservation_2018, seibold_m-matrices_2006}. We therefore repeat the stability analysis for grids in which there is a minimum distance between grid points. This is implemented using the `randomness' parameter r \eqref{eq:randomness}. The results are given in table \ref{tab:1DGridSensitivity}.

As expected, the first order upwind scheme is consistently stable \eqref{eq:1DCFLCondition}. The second order upwind and the second order MUSCL scheme are also stable for all the checked grids. The first order and third-order MUSCL schemes always have at least one eigenvalue in the right half plane, and are therefore not included in the table. The fourth order MUSCL scheme is more frequently unstable. As the amount of grid points is increased, the scheme is more likely to be unstable. As the grid is more regular (\(r\) decreases), the stability of the scheme improves. In the case of \(r = \frac{8\Delta x}{20}\), all schemes were stable for all 100 grids that were tested. 

For the case \(N = 100\) and \(r = \frac{\Delta x}{2}\), actual long-time simulations were performed with a very small time step (CFL = 0.05) using the exact same grids that were used to generate the tables. We saw instabilities in the solutions for the same grids for which we found unstable eigenvalues, thereby confirming our results. We also found that the WENO scheme was unstable for four grids. The meshless WENO scheme is therefore more unstable than the two newly proposed MUSCL schemes.
\begin{remark}
	The stability of meshless schemes strongly depends on the weight function that is used in the underlying least squares problem \eqref{eq:WeightFunction}. See \cite{tey_moving_2021} for an overview of existing weight functions.  
\end{remark} 

\begin{table}
	\begin{subtable}{1\textwidth}
		\sisetup{table-format=-1.2}   
		\centering
		\begin{tabular}{|c||c|c|c|c|} 
			\hline
			$N$ & Upwind 1 & Upwind 2 & MUSCL 2 & MUSCL 4 \\ \hline \hline
			100 & 0.0 & 0.0 & 0.0 & 2.0 \\
			200 & 0.0 & 0.0 & 0.0 & 5.0 \\
			300 & 0.0 & 0.0 & 0.0 & 5.0 \\
			400 & 0.0 & 0.0 & 0.0 & 5.0 \\
			\hline
		\end{tabular}
		\caption{$r = \frac{\Delta x}{2}$}\label{tab:sub_first}
	\end{subtable}
	
	\bigskip
	\begin{subtable}{1\textwidth}
		\sisetup{table-format=4.0} 
		\centering
		\begin{tabular}{|c||c|c|c|c|} 
			\hline
			$N$ & Upwind 1 & Upwind 2 & MUSCL 2 & MUSCL 4 \\ \hline \hline
			100 & 0.0  & 0.0 & 0.0 & 0.0 \\
			200 & 0.0  & 0.0 & 0.0 & 1.0 \\
			300 & 0.0  & 0.0 & 0.0 & 2.0 \\
			400 & 0.0  & 0.0 & 1.0 & 3.0 \\
			\hline
		\end{tabular}
		\caption{$r = \frac{9 \Delta x}{20}$}\label{tab:sub_second}
	\end{subtable}
	\caption{Percentage of the number of times (out of 100) that the scheme was unstable. The first column is the amount of grid points $N$.} 
	\label{tab:1DGridSensitivity}
\end{table}

\subsubsection{Two space dimensions}
\begin{figure}[h!]
	\centering
	\includegraphics[width=0.8\linewidth]{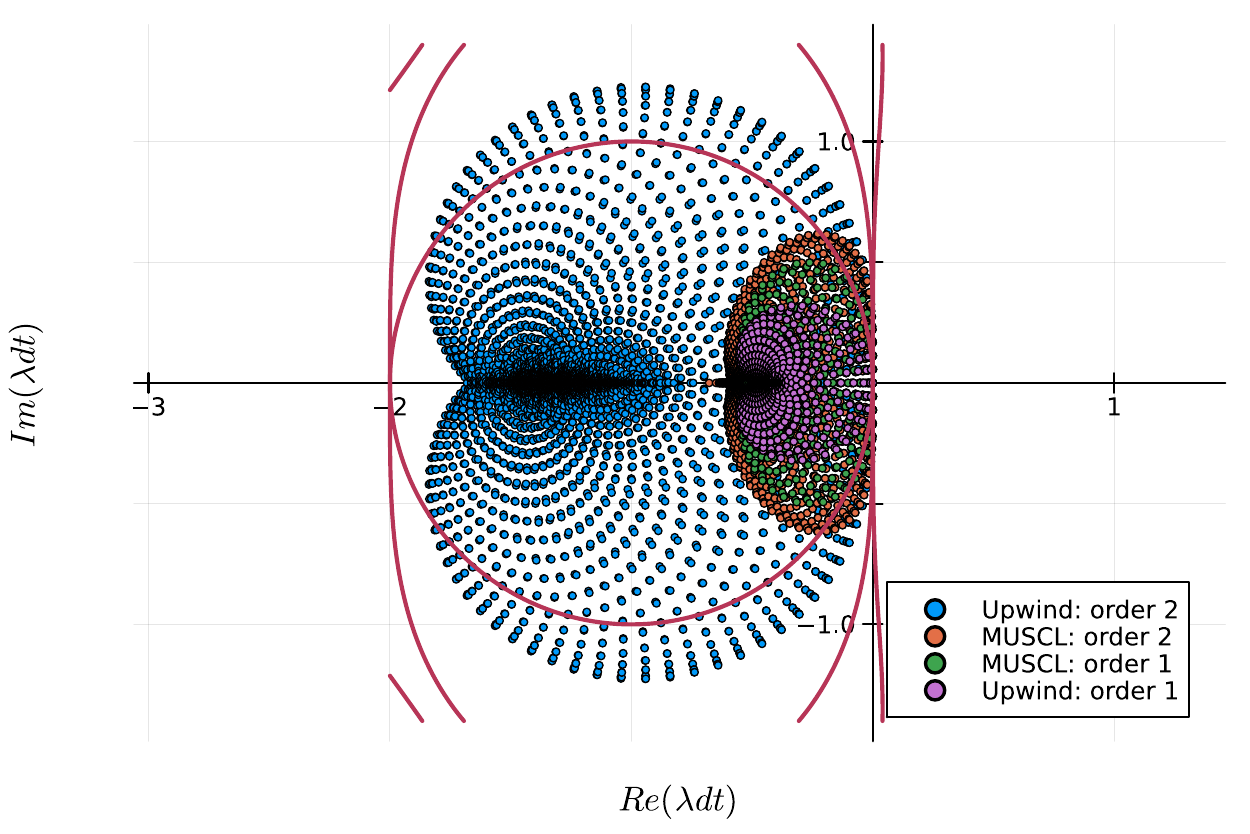}
	\caption{Spectra of several meshless schemes with CFL number $1/2$. The lines in the graph are the stability regions for the first, second, third and fourth order explicit Runge-Kutta method.}
	\label{fig:2DEigenvalues}
\end{figure}

We repeat the analysis in 2D. First, we consider a irregular grid $\Omega = \left[ -5, 5\right]^2$ with \(70^2\) grid points and plot the spectra of the semi-discretised ODEs. In figure \ref{fig:2DEigenvalues}, the spectra of several meshless schemes are plotted. We note that the randomness parameter used to generate the grids is set to the maximum value of $r = \frac{\Delta x}{2}$. The first and second order MUSCL and upwind schemes have stable eigenvalues. As in the 1D case, the second order upwind scheme has the largest spectral radius, and the second order MUSCL scheme has the most eigenvalues that extend along the imaginary axis. We repeat this analysis for 100 grids and log if unstable eigenvalues are found. Eigenvalues are considered unstable if the imaginary part is larger than \(1\times10^{-13}\). We crosscheck these results by performing a long-time simulation (\(t_{max} = 30\sqrt{2}, \text{CFL} = 0.1\)) with a shock initial condition. This allows us to also check the stability of the WENO scheme. For all grids, the first-and second order MUSCL schemes, the first order upwind scheme from section \ref{section:FirstOrder2D}, and the second order upwind scheme were always stable. Note that in one dimension, the first order MUSCL scheme was always unstable. The WENO scheme was only stable in \(72\%\) of the checked grids.

\subsection{Dirichlet boundaries}
\label{section:DBoundary}
In this section, we examine the behavior of the numerical methods in the presence of a Dirichlet boundary. We consider a one-dimensional non-periodic domain, \(\Omega = \left[-5, 5\right]\), with a shock located at \(x = -2\). On the left boundary, we impose a Dirichlet condition \(u(t, -5) = 0.5\). Consequently, near the boundary, the central stencil used in the MUSCL schemes becomes biased toward one side. For the WENO method, one-sided stencils with an insufficient number of points are deactivated near the boundary, forcing the method to rely on the central and remaining one-sided stencils. These slight modifications near the boundary do not pose any issue as long as the solution remains smooth or an upwind stencil is available.

\begin{figure}[h!]
	\centering
	\includegraphics[width=0.7\linewidth]{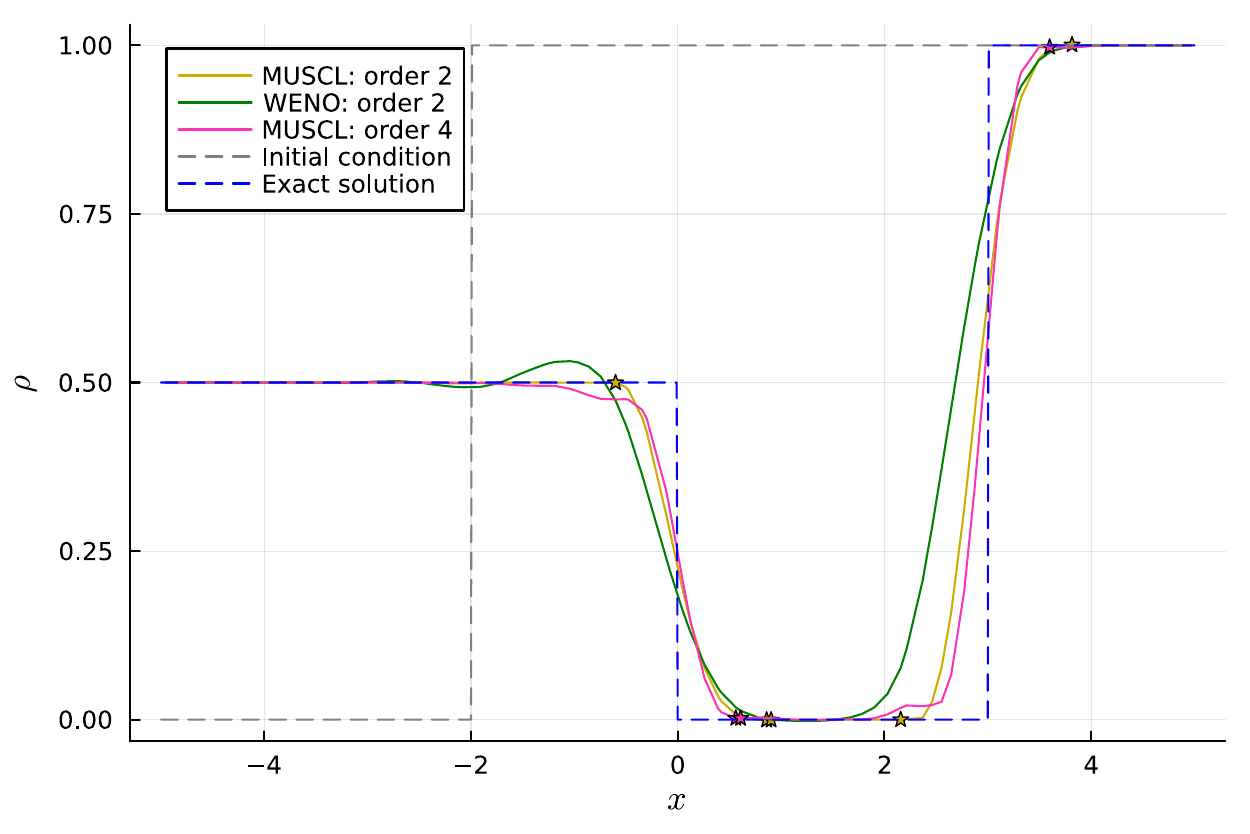}
	\caption{A simulation with a Dirichlet boundary on the left. The simulation is performed with 100 grid points and CFL number \(1/3\).}
	\label{fig:DBoundary}
\end{figure}

We perform a simulation up to time \(t = 5\) using 100 grid points and a CFL number of \(1/3\). We use the same combination of algorithms as in section \ref{section:Convergence}. The results are shown in figure \ref{fig:DBoundary}. The WENO scheme produces a smooth solution profile with slight overshoot, while the second- and fourth-order MUSCL schemes yield sharper profiles. As expected, no oscillations appear at the left boundary where the Dirichlet condition is applied.

\subsection{Conservation}
\label{section:Conservation}

\begin{figure}
	\centering
	\includegraphics[width=0.5\linewidth]{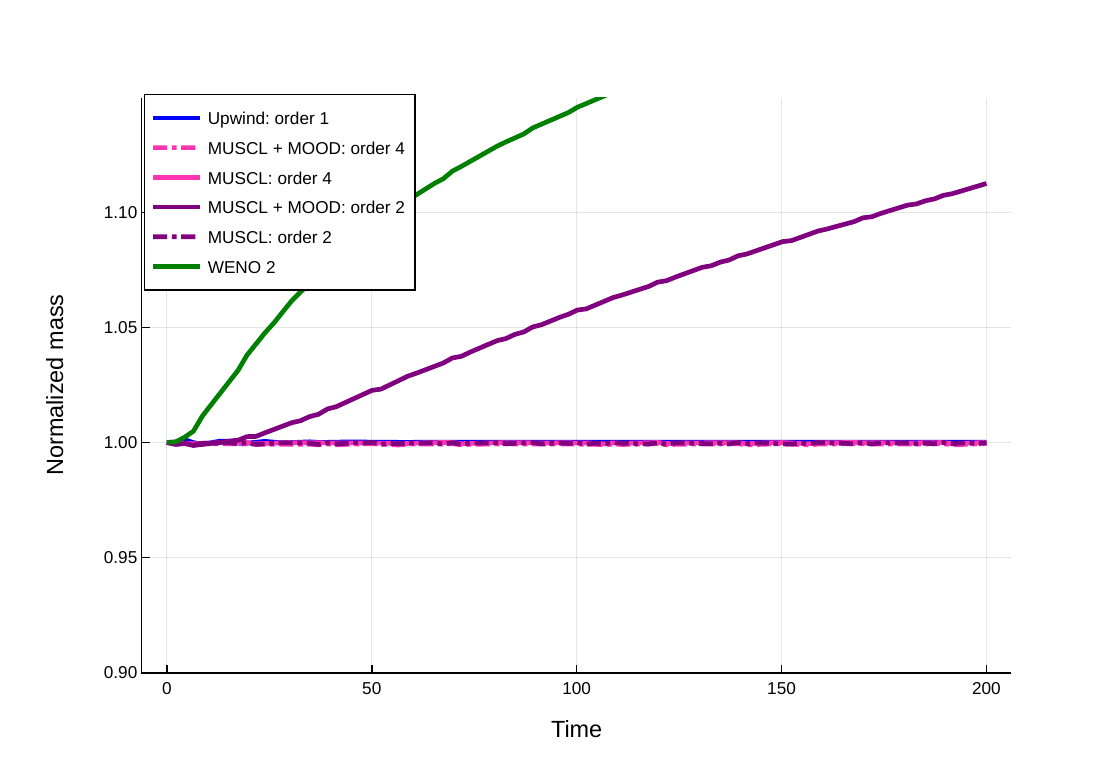}
	\caption{Normalised mass as a function of time for several meshless schemes.}
	\label{fig:1DMassLoss}
\end{figure}

It was already mentioned that the MUSCL schemes cannot be conservative due to the lack of flux conservation and geometric conservation, see section \ref{section:conservationRemark}. In this section, we quantify the lack of conservation by performing a long-time simulation and tracking the mass in the system. We consider a random fixed irregular grid of $N = 100$ grids points on the domain $\Omega \in \left[ -5, 5\right]$ with a smooth Gaussian condition \begin{align} u(0, x) = \exp(-x^2). \end{align} We consider the WENO scheme, the first order upwind scheme, and several MUSCL schemes. All schemes are paired with a third-order Runge-Kutta time integration routine. We do simulations up to $t = 200$ with the CFL number equal to $\frac{1}{4}$ such that all schemes are stable. The total mass is computed using a first order quadrature rule. The normalised mass as a function of time is plotted in figure \ref{fig:1DMassLoss}. 

The WENO scheme and second order MUSCL-MOOD scheme, add a significant amount of mass. Note that the second order MUSCL scheme without MOOD, and the first order upwind scheme, preserve the mass surprisingly well, but their combination in the second order MUSCL scheme with MOOD, increases the mass over time. From the solution plots, it can be seen that the solution due to the second order MUSCL scheme becomes negative, yet the total mass in the system remains constant. The increase of mass in the second order MUSCL 2 scheme is what makes the scheme positive. Interestingly, the fourth order MUSCL-MOOD scheme does not exhibit this behavior, and maintains the mass in the system. Note that all schemes here were stable; no oscillations were seen in the solution. The increase or decrease of mass is thus purely due to the scheme and is not due to a lack of linear stability, which could also blow up the mass. 

\subsection{Efficiency}
\label{section:Efficiency}

\subsubsection{One space dimension}
\begin{figure}[h!]
	\centering
	\begin{subfigure}[b]{0.49\textwidth}
		\centering
		\includegraphics[width=\linewidth]{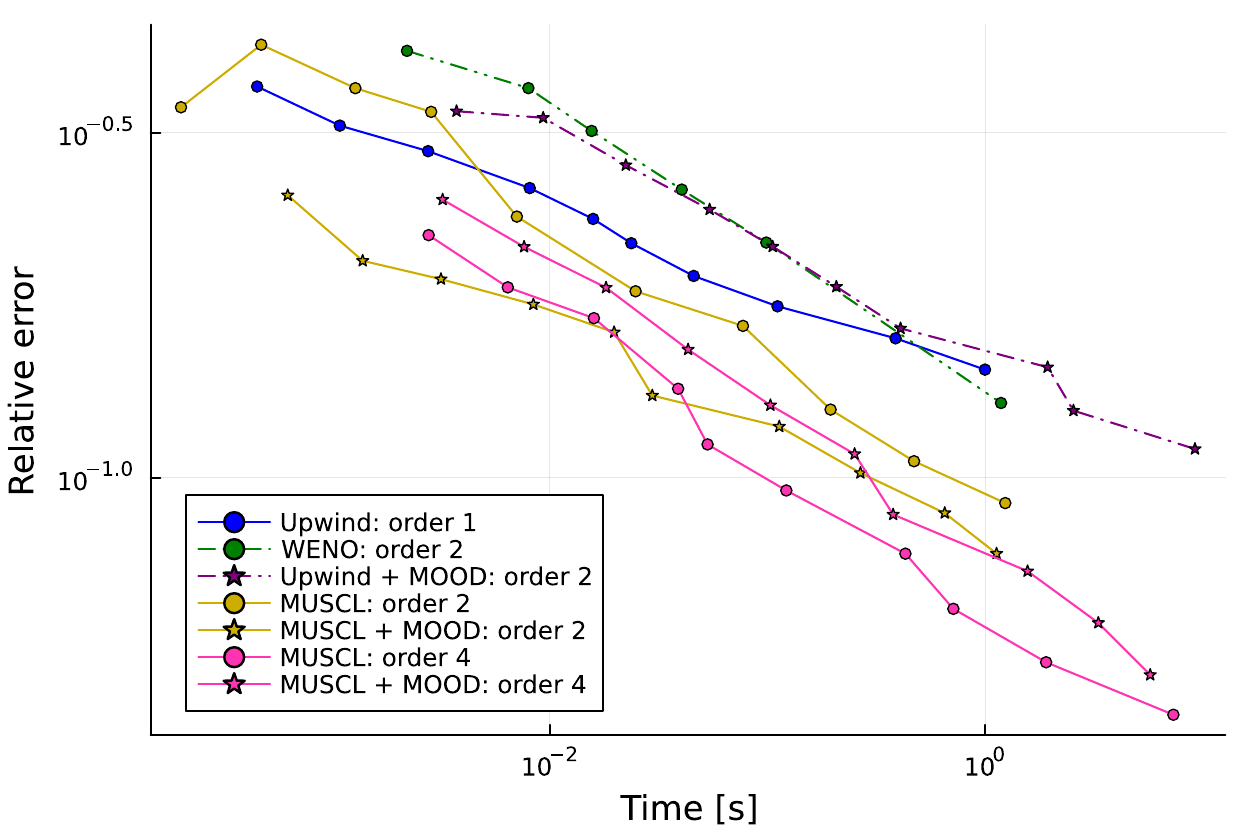}
		\label{fig:1DEfficiencyShock}
	\end{subfigure}
	\begin{subfigure}[b]{0.49\textwidth}
		\centering
		\includegraphics[width=\linewidth]{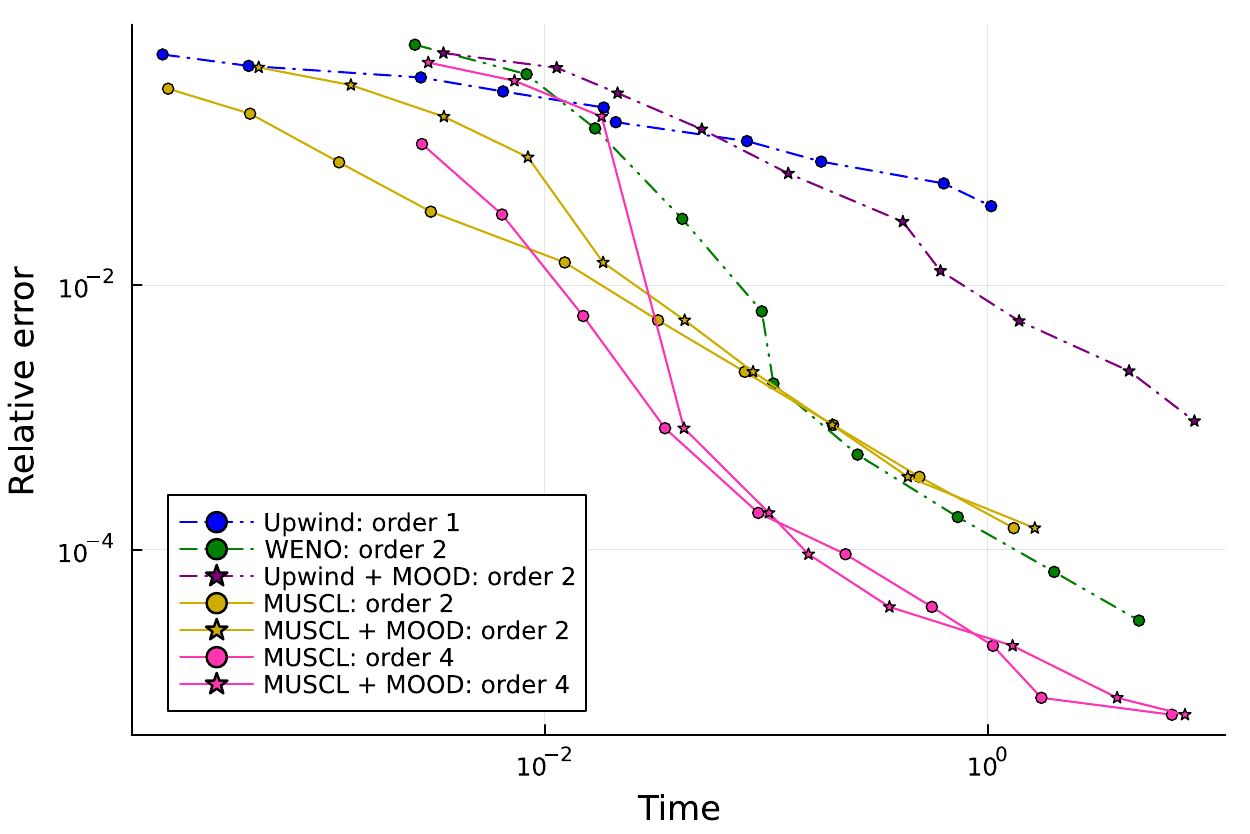}
		\label{fig:1DEfficiencySmooth}
	\end{subfigure}
	\caption{Efficiency of numerical schemes (see table \ref{tab:1DEfficiencyTable}) for a shock (left) and smooth (right) initial condition. Simulations are performed for grid sizes 30, 46, 72, 111, 171, 264, 407, 629, 971 and 1500. Some noise on the timings is still visible.}
	\label{fig:1DEfficiency}
\end{figure}

To conclude the analysis of the new MUSCL schemes, we report on the algorithm's efficiency: the computational time required for the algorithm to obtain the solution for a given error threshold. We again consider a periodic domain $\Omega = \left[ -5, 5\right]$ with a fixed irregular grid. We solve the linear advection equation up to $t = 7.5$, for a smooth initial condition and a shock initial condition \begin{align}
	u_1(0, x) = \exp \left(-x^2\right), \; \; u_2(0, x) = \begin{cases}
		1 & \text{if } x > 0 \\
		0 & \text{else. }
	\end{cases}
\end{align}
The simulations are performed with the first order and second order upwind methods, the second order WENO, and the second-and fourth order MUSCL methods. To obtain the most realistic results, we match the order of the time integrator with the order of the spatial discretisation method. For all higher-order schemes except the WENO method, we use a MOOD procedure with a fallback to the first order upwind method. For each algorithm, we use the largest possible time step for which the scheme is stable. We summarize the combinations of schemes and the time step in \ref{tab:1DEfficiencyTable}. Finally, time and error results are averaged over ten random grids to reduce the noise in the plot. The efficiency of the schemes is plotted in figure \ref{fig:1DEfficiency}.

\begin{table}
	\centering
	\begin{adjustbox}{width=\columnwidth,center}
		\begin{tabular}{|c||c|c|c|c|c|c|c|} 
			\hline
			Time & Forward Euler & Ralston RK2 & Ralston RK2 & Ralston RK2 & RK4 & Ralston RK2 & RK4 \\ \hline 
			Space & Upwind 1 & Upwind 2 & WENO 2 & MUSCL 2 & MUSCL 4 & MUSCL 2 & MUSCL 4\\ \hline 
			MOOD & \xmark & \cmark & \xmark & \cmark & \cmark & \xmark & \xmark \\ \hline
			CFL & 0.99 & 0.3 & 0.7 & 0.75 & 0.7 & 0.75 & 0.7 \\ \hline
		\end{tabular}
	\end{adjustbox}
	\caption{Combinations of algorithms used in the efficiency test. The CFL number is relative to the maximal allowed time step of the first order Euler and Upwind combination: $dt = \text{CFL} dt_{euler}$, where $dt_{euler}$ is the time step for which the equality in \eqref{eq:1DCFLCondition} is satisfied.}
	\label{tab:1DEfficiencyTable}
\end{table}

The reference first order upwind scheme, and the second order upwind scheme prove to be inefficient. The WENO and second order MUSCL scheme with MOOD achieve similar results, the former being better for the smooth initial condition, the latter being better for the shock initial condition. It should be noted that the WENO scheme was unstable for some of the test grids over which the results were averaged for $N_x = 1500, 971 \text{and } 407$. These results were removed from the plot, but it nonetheless illustrates the unreliable nature of the WENO scheme. The best results are obtained with the fourth order MUSCL scheme with MOOD. This scheme is the most efficient already for very low error requirements. From the graph it can be seen that the additional computational cost of the MOOD procedure is limited.

\subsubsection{Two space dimensions}
\begin{figure}[h!]
	\centering
	\begin{subfigure}[b]{0.49\textwidth}
		\centering
		\includegraphics[width=\linewidth]{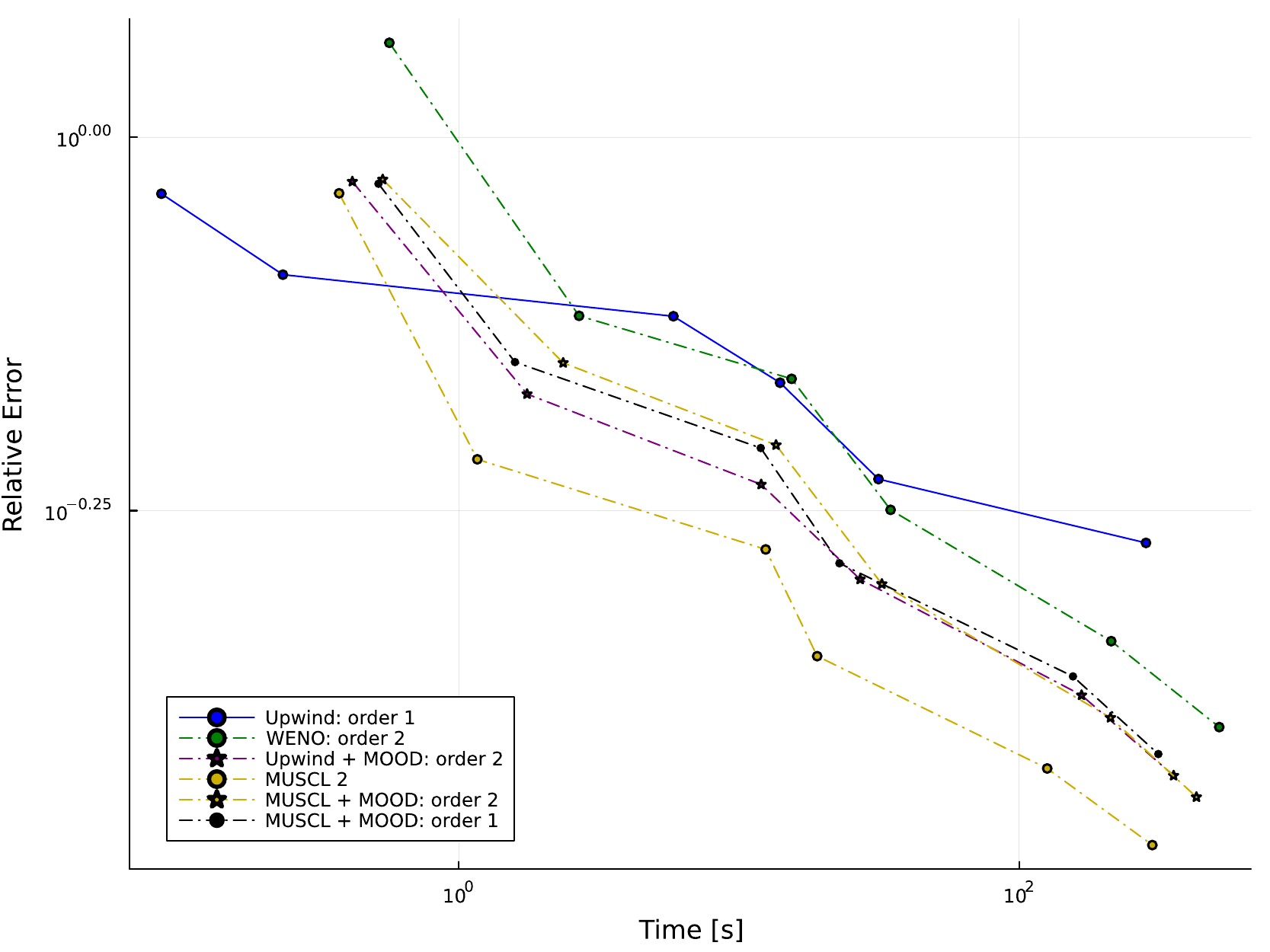}
		\label{fig:2DEfficiencyShock}
	\end{subfigure}
	\begin{subfigure}[b]{0.49\textwidth}
		\centering
		\includegraphics[width=\linewidth]{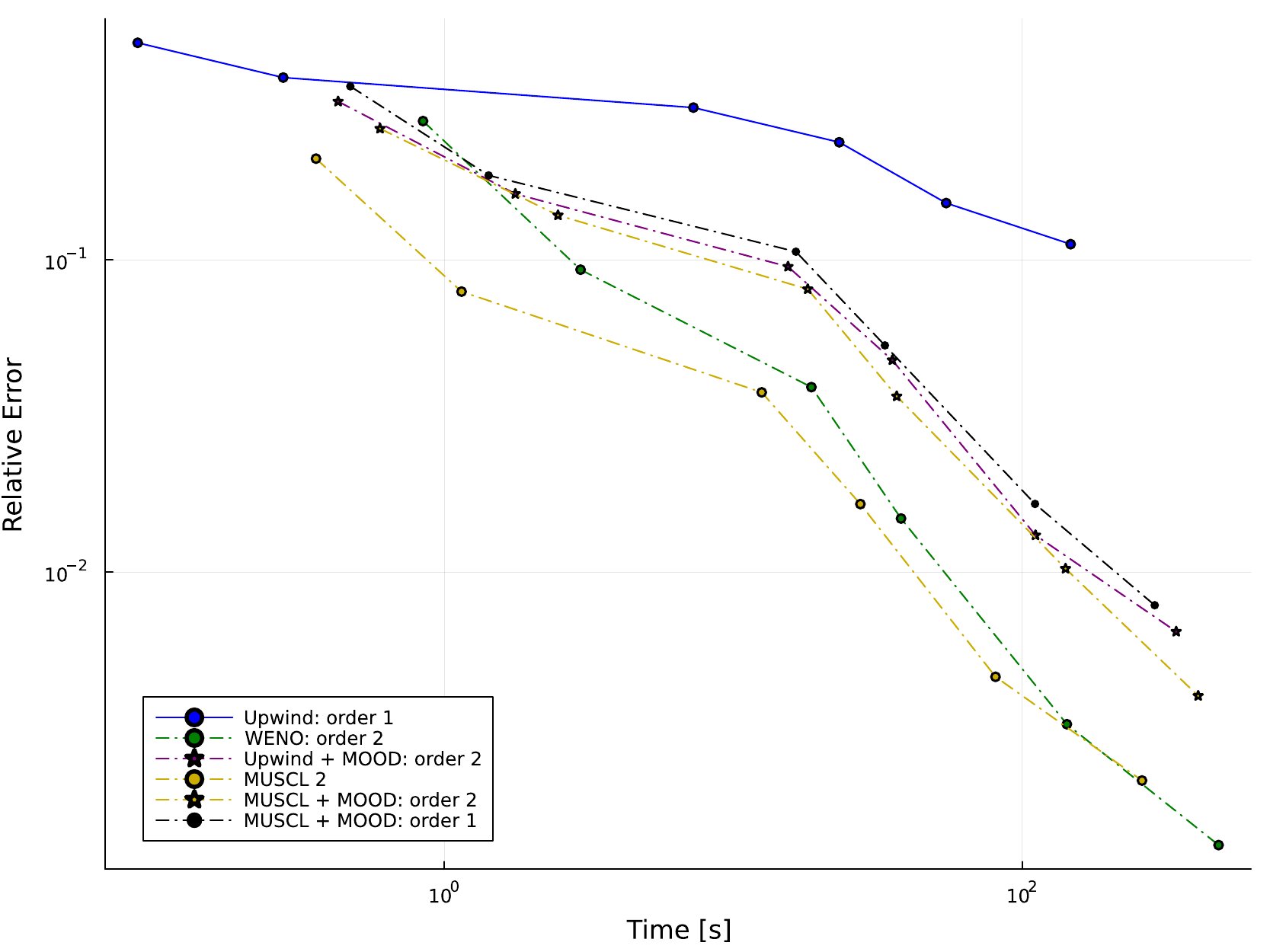}
		\label{fig:2DEfficiencySmooth}
	\end{subfigure}
	\caption{Efficiency of numerical schemes for a shock (left) and smooth (right) initial condition. Simulations are performed for grid sizes \(30^2\), \(50^2\), \(70^2\), \(100^2\), \(175^2\) and \(250^2\).}
	\label{fig:2DEfficiency}
\end{figure}

We repeat the analysis in 2D. We consider a fixed irregular grids for the periodic domain $\Omega = \left[ -5, 5\right]^2$. We solve the linear advection equation up to time \(t = 1.0\) for a smooth and shock initial condition \eqref{eq:2DSmoothShockInit}. All schemes, except for the first order method, are paired with a third-order Runge-Kutta time integration scheme. The CFL number is fixed to 0.5. The results are plotted in figure \ref{fig:2DEfficiency}. 
For the smooth initial condition, the second order WENO scheme and the second order MUSCL scheme perform equally well. The second order MUSCL scheme with MOOD does not achieve the same performance as the `plain' second order MUSCL scheme for two reasons. Firstly, checking the MOOD criteria for each point adds 5-10\% computation time. As a result, the efficiency line of the second order MUSCL scheme with MOOD is shifted to the right with respect to the same algorithm without MOOD. Secondly, even for the smooth initial condition, there are still a very small amount of MOOD events occurring at the peak of the Gaussian throughout time integration, which results in an increase of the error. The first order upwind scheme, second order upwind scheme, and first order MUSCL scheme do not yield as good results. For the shock initial condition, the best results are obtained with the first order MUSCL scheme with MOOD, second order MUSCL scheme with MOOD, and the second order upwind scheme with MOOD. 

\section{Conclusion}
\label{section:Conclusion}
In this paper, we discussed the order, efficiency, stability and positivity of several meshless schemes for linear scalar hyperbolic equations. In the comparison, we consider existing upwind schemes and WENO schemes, and a new class of MUSCL-like schemes. The new MUSCL-like meshless scheme uses a central stencil and can achieve arbitrary high orders. The stability of the scheme is guaranteed by an upwind reconstruction to the midpoints of the stencil. The new MUSCL schemes are also efficient due to the reuse of the GFDM solution in the reconstruction. The new MUSCL schemes are combined with a Multi-dimensional Optimal Order Detection (MOOD) procedure to obtain a scheme that does not yield spurious oscillations at discontinuities. In one spatial dimension, our fourth order MUSCL scheme outperforms existing WENO and upwind schemes in terms of stability and accuracy. In two spatial dimensions, our MUSCL scheme achieves similar accuracy to an existing WENO scheme but is significantly more reliable. In future work, these schemes will be extended to non-linear (systems) of conservation laws with moving meshes. In addition, we intend to use the meshless MUSCL scheme to solve kinetic equations on moving irregular grids.

\textbf{Declarations.} The authors have received funding from the European Union's Framework Program for Research and Innovation HORIZON-MSCA-2021-DN-01 under the Marie Sklodowska-Curie Grant Agreement Project 101072546 – DATAHYKING.

\textbf{Data Availability Statement.} All code and results are available at \cite{githubRepo}.

\newpage
\bibliographystyle{acm}
\bibliography{Klaas, manualBib}

\end{document}